\documentclass{amsart}[11pt]
\usepackage{amssymb}
\usepackage[dvips]{color}
\usepackage{graphicx}
\usepackage{epsfig}

\newcommand\strict{\operatorname{strict}}
\DeclareMathOperator{\weak}{{\mathrm{weak}}}

\def\SG{{\mathfrak S}}
\def\A{{\mathcal A}}
\def\Comp{{\rm Comp}}
\def\Z{\mathbb Z}
\def\QSym{{\rm QSym}}
\def\a{{\mathbf a}}
\def\b{{\mathbf b}}
\def\c{{\mathbf c}}

\def\deg{{\rm deg}}

\def\J{\mathcal J}

\def\C{\mathcal C}

\def\P{\mathbb P}
\def\N{\mathbb N}
\def\T{\mathbb T}
\def\Z{\mathbb Z}
\def\wt{\mathrm{wt}}

\newtheorem{theorem}{Theorem}[section]
\newtheorem{lemma}[theorem]{Lemma}
\newtheorem{proposition}[theorem]{Proposition}
\newtheorem{corollary}[theorem]{Corollary}
\newtheorem{conjecture}[theorem]{Conjecture}

\theoremstyle{definition}
\newtheorem{example}[theorem]{Example}
\newtheorem{problem}[theorem]{Problem}
\newtheorem{definition}[theorem]{Definition}
\newtheorem{question}[theorem]{Question}
\theoremstyle{remark}

\newtheorem{remark}[theorem]{Remark}

\def\weak{\rm weak}
\def\strict{\rm strict}
\def\p{\mathbf{p}}
\def\N{\mathbb N}
\def\Z{\mathbb Z}
\begin{document}

\title{$P$-partition products and \\ fundamental quasi-symmetric function positivity}
\author{Thomas Lam and Pavlo Pylyavskyy}
\thanks{T.L. was partially supported by NSF DMS-0600677.}
\thanks{We thank R.~Stanley and A.~Postnikov for interesting
discussions.}

\address{T.L.: Department of Mathematics, Harvard, Cambridge, MA, 02138}
\email{tfylam@math.harvard.edu}

\address{P.P.: Department of Mathematics, M.I.T., Cambridge, MA, 02139}
\email{pasha@mit.edu}

\date{September 6, 2006}

\begin{abstract}
We show that certain differences of products $$K_{Q \wedge R,\theta}
\, K_{Q\vee R,\theta} - K_{Q,\theta} \, K_{R,\theta}$$ of
$P$-partition generating functions are positive in the basis of
fundamental quasi-symmetric functions $L_\alpha$.  This result
interpolates between recent Schur positivity and monomial positivity
results of the same flavor. We study the case of chains in detail, introducing
certain ``cell transfer'' operations on compositions and an interesting related 
``$L$-positivity'' poset.
We introduce and study quasi-symmetric functions called {\it wave Schur functions} and
use them to establish, in the case of chains, that $K_{Q \wedge R,\theta}
\, K_{Q\vee R,\theta} - K_{Q,\theta} \, K_{R,\theta}$ is itself equal to a single generating function
$K_{P,\theta}$ for a labeled poset $(P,\theta)$.  In the course of our investigations we
establish some factorization properties of the ring $\QSym$ of quasisymmetric functions.
\end{abstract}

\maketitle

\markboth{Lam and Pylyavskyy}{$L$-positivity} \maketitle

\section{Introduction} \label{sec:intro}
The present article studies quasi-symmetric functions $f \in \QSym$
which are non-negative linear combinations of fundamental
quasi-symmetric functions.  It sits between joint work with
Postnikov~\cite{LPP} involving symmetric functions and Schur
positivity and our more poset-theoretic work~\cite{LP} involving
monomial positivity.

\medskip
The Schur functions $s_\lambda$ (see~\cite{EC2}) are important symmetric functions
which occur in combinatorics, representation theory and geometry.
They exhibit a multitude of remarkable properties, and here we
highlight three non-trivial {\it positivity} properties.
\begin{itemize}
\item[(A)]
The product $s_\lambda \, s_\mu$ of two Schur functions is
Schur-positive (Littlewood-Richardson rule I).
\item[(B)]
A skew Schur function $s_{\lambda/\mu}$ is Schur-positive
(Littlewood-Richardson rule II).
\item[(C)]
The difference of products $s_{\max(\lambda/\mu,\nu/\rho)}\,
s_{\min(\lambda/\mu,\nu/\rho)} - s_{\lambda/\mu}\,
s_{\nu/\rho}$ is Schur positive
(Lam-Postnikov-Pylyavskyy~\cite{LPP}).  Here $\max$ and $\min$ are
taken coordinate-wise.
\end{itemize}

The aim of this article is to replace symmetric functions with
quasi-symmetric functions and study analogous positivity properties
for the fundamental quasi-symmetric functions $L_\alpha$.  Say that
a quasi-symmetric function $f$ is {\it $L$-positive} if it is a
non-negative linear combination of fundamental quasi-symmetric
functions.  We take the point of view that the quasi-symmetric
analogues of (A) and (B) are:
\begin{itemize}
\item[(A*)]
The product $L_\alpha L_\beta$ of fundamental quasi-symmetric
functions is $L$-positive (shuffle product).
\item[(B*)]
For any poset $P$ and labeling $\theta: P \to \P$ the generating
function $K_{P,\theta}$ of $P$-partitions is $L$-positive (Stanley's
$P$-partition theory~\cite{Sta}).
\end{itemize}
Thus the functions $K_{P,\theta}$ will replace the skew Schur
functions $s_{\lambda/\mu}$.  Our main result
(Theorem~\ref{thm:main}), which is the analogue of (C), says that
the difference
\begin{equation} \label{eq:dif} K_{Q \wedge R,\theta} \, K_{Q\vee
R,\theta} - K_{Q,\theta} \, K_{R,\theta}\end{equation} is
$L$-positive, where $Q$ and $R$ are two convex subsets of a labeled
poset $(P,\theta)$ and $\wedge$ and $\vee$ are the {\it cell
transfer operations} introduced in~\cite{LP}.

\medskip
In~\cite{LP}, the same difference (\ref{eq:dif}) is shown to be
monomial-positive for a larger class of posets called $\T$-labeled
posets.  Since Schur-positivity implies $L$-positivity which in turn
implies monomial-positivity, our current result sits between the two
results of~\cite{LP} and~\cite{LPP}; with each restriction to the
class of (labeled) posets a stronger form of positivity holds.  We
summarize the relationship between this article and the two
previous works~\cite{LP,LPP} in a table.
\begin{center}
\begin{tabular}{|c||c|c|c|}
\hline Paper & Cell Transfer~\cite{LP} & This& Schur
positivity~\cite{LPP}
\\
\hline Ring & $\Z[[x_1,x_2,\ldots]]$ & $\QSym$ & ${\rm Sym}$\\
\hline Basis & $x^\alpha$ &
$L_\alpha$ & $s_\lambda$ \\
\hline Skew fcns. & $K_{P,O}$ & $K_{P,\theta}$ &
$s_{\lambda/\mu}$ \\
\hline Posets & $\T$-labeled posets $(P,O)$ & Stanley's $(P,\theta)$
& Young
diagrams $\lambda/\mu$ \\
\hline
\end{tabular}
\end{center}

In each of the three cases, the structure constants in the ``basis'' are non-negative,
and the ``skew functions'' lie in the
``ring'' and are non-negative when written in terms of the
``basis''.  In all three cases, the difference of products of ``skew
functions'' arising from the cell transfer operation on the
``posets'' is positive in the corresponding ``basis''.  The Schur functions have
an interpretation as irreducible characters of the symmetric group while
the fundamental quasi-symmetric functions have an interpretation as
irreducible characters of the 0-Hecke algebra.  It would be interesting to
give a representation theoretic explanation of our results (and in particular 
of the cell transfer operations).

\medskip

We study the difference (\ref{eq:dif}) in detail for the special
case that $Q$ and $R$ are convex subsets of a chain, in which case
all the four functions in this difference are themselves fundamental
quasi-symmetric functions $L_\alpha$.  We introduce ``cell
transfer'' operations on compositions, also denoted $\vee$ and
$\wedge$, such that the difference $L_{\alpha \wedge \beta} \,
L_{\alpha \vee \beta} - L_\alpha \, L_\beta$ is $L$-positive.  We
further conjecture (Conjecture~\ref{conj:max}) that a product
$L_\alpha\, L_\beta$ is maximal in ``$L$-positivity order'' if and
only if the pair $\{\alpha, \beta\}$ is stable under cell transfer.
As part of this investigation, we show that each $L_\alpha$ is
irreducible.

Next, we ask when the difference
(\ref{eq:dif}) is itself equal to $K_{P,\theta}$ for some labeled
poset $(P,\theta)$.  We show that this is always the case for the
differences $L_{\alpha \wedge \beta} \, L_{\alpha \vee \beta} -
L_\alpha \, L_\beta$ by introducing generating functions we call
{\it wave Schur functions}, which appear to be interesting in 
their own right.

Wave Schur functions are generating functions of certain Young
tableaux, where the weakly and strictly increasing conditions are
altered in a particular alternating pattern.  We call these tableaux
``{\it wave $\p$-tableaux}'' where $\p$ indicates how the increasing
conditions have been modified.  We show that wave Schur functions
are $L$-positive and that they
satisfy a Jacobi-Trudi style determinantal formula
(Theorem~\ref{jtl}), with the fundamental quasi-symmetric functions
replacing the homogeneous symmetric functions.  The difference
$L_{\alpha \wedge \beta} \, L_{\alpha \vee \beta} - L_\alpha \,
L_\beta$ is equal to an appropriate wave Schur function for a
two-row skew shape.

\medskip
In the final sections of the paper, we comment on whether our
results can be expanded to a larger class of generating functions of
the form $K_{P,O}$ for a $\T$-labeled poset $(P,O)$.

\section{Quasi-symmetric functions}
\label{sec:sym} We refer to~\cite{EC2} for more details of the
material in this section.

\subsection{Monomial and fundamental quasi-symmetric functions}
Let $n$ be a positive integer.  A {\it composition} of $n$ is a
sequence $\alpha = (\alpha_1,\alpha_2,\ldots,\alpha_k)$ of positive
integers such that $\alpha_1 + \alpha_2 + \cdots + \alpha_k = n$. We
write $|\alpha| = n$.  Denote the set of compositions of $n$ by
$\Comp(n)$.  Associated to a composition $\alpha =
(\alpha_1,\alpha_2,\ldots,\alpha_k)$ of $n$ is a subset $D(\alpha) =
\{\alpha_1,\alpha_1+\alpha_2,\ldots,\alpha_1+\alpha_2+\dots +
\alpha_{k-1}$ of $[n-1]$.  The map $\alpha \mapsto D(\alpha)$ is a
bijection between compositions of $n$ and subsets of $[n-1] =
\{1,2,\ldots,n-1\}$.  We will denote the inverse map by $\C:
2^{[n-1]} \to \Comp(n)$ so that $\C(D(\alpha))= \alpha$.

%If in addition $\alpha_1 \geq \alpha_2 \geq \cdots \geq \alpha_k$
%then we say that $\alpha$ is a {\it partition} of $n$.  If $\lambda$
%is a partition then $\lambda'$ denotes the {\it conjugate}
%partition.  Let $l(\ll)$ denote the number of (non-zero) parts of
%$\lambda$.

A formal power series $f = f(x) \in \Z[[x_1,x_2,\ldots]]$ with
bounded degree is called {\it quasi-symmetric} if for any
$a_1,a_2,\ldots,a_k \in \P$ we have
\[
\left[x_{i_1}^{a_1}\cdots x_{i_k}^{a_k} \right]f =
\left[x_{j_1}^{a_1}\cdots x_{j_k}^{a_k} \right]f
\]
whenever $i_1 < \cdots < i_k$ and $j_1 <\cdots < j_k$.  Here
$[x^\alpha]f$ denotes the coefficient of $x^\alpha$ in $f$. Denote
by $\QSym \subset  \Z[[x_1,x_2,\ldots]]$ the ring of quasi-symmetric
functions.

Let $\alpha$ be a composition. Then the {\it monomial
quasi-symmetric function} $M_\alpha$ is given by
\[
M_\alpha = \sum_{i_1 < \cdots <i_k} x_{i_1}^{\alpha_k} \cdots
x_{i_k}^{\alpha_k}.
\]
The {\it fundamental quasi-symmetric function} $L_\alpha$ is given
by
\[
L_\alpha = \sum_{D(\beta) \subset D(\alpha)} M_\beta,
\]
where the summation is over compositions $\beta$ satisfying $|\beta| = |\alpha|$.
The set of fundamental quasi-symmetric functions (resp. monomial
quasi-symmetric functions) form a basis of $\QSym$.  We say that a
quasi-symmetric function $f \in \QSym$ is {\it $L$-positive} (resp.
{\it $M$-positive}) if it is a non-negative linear combination of
fundamental quasi-symmetric functions (resp. monomial
quasi-symmetric functions).  Note that $L$-positivity implies
$M$-positivity.

Two fundamental quasi-symmetric functions $L_\alpha$ and $L_\beta$
multiply according to the {\it shuffle product}.  Let $u = u_1 u_2
\cdots u_k$ and $v = v_1 v_2 \cdots v_l$ be two words.  Then a word
$w = w_1 w_2 \cdots w_{k+l}$ is a shuffle of $u$ and $v$ if there
exist disjoint subsets $A, B \subset [k+l]$ such that $A =
\{a_1,a_2,\ldots,a_k\}$, $B = \{b_1,b_2,\ldots,b_l\}$, $w_{a_i} =
u_i$ for all $1 \leq i \leq k$, $w_{b_i}= v_i$ for all $1 \leq i
\leq l$ and $A \cup B = [k+l]$.  We denote the set of shuffles of
$u$ and $v$ by $u \odot v$.

For a composition $\alpha$ with $|\alpha| = n$ let $w(\alpha) = w =
w_1 w_2 \cdots w_n$ denote any word with {\it descent set} $D(w) =
\{i : w_i > w_{i+1}\} \subset [n-1]$ equal to $D(\alpha)$.  Suppose
$w(\alpha)$ and $w(\beta)$ are chosen to have disjoint letters.
Then
$$
L_\alpha \, L_\beta = \sum_{u \in w(\alpha) \odot w(\beta)}
L_{\C(u)},
$$
where $\C(u)$ is by definition the composition $\C(D(u))$ associated
to $u$.

\subsection{Two involutions on $\QSym$}
If $D \subset[n-1]$ we let $\bar D = \{i \in[n-1] \mid i \notin D\}$
denote its complement.  For a composition $\alpha$, define $\bar
\alpha = \C(\overline{D(\alpha)})$.  Let $\omega$ denote the linear
endomorphism of $\QSym$ given by $\omega(L_\alpha) = L_{\bar
\alpha}$.

Let $\alpha^*$ denote $\alpha$ read backwards: $\alpha^* =
(\alpha_k, \ldots, \alpha_1)$. Let $\nu$ be the linear endomorphism
of $\QSym$ which sends $L_{\alpha} \mapsto L_{\alpha^*}$.

\begin{proposition}
\label{prop:inv} The maps $\omega$ and $\nu$ are algebra involutions
of $\QSym$, and we have $\nu(M_{\alpha})= M_{\alpha^*}$.
\end{proposition}

\begin{proof}
We will check the first statement for $\nu$; the proof for $\omega$
is similar.  Let $w =w_1w_2\cdots w_r \in S_r$ be a permutation with
descent set $D(w) = D(\alpha)$. Then $w^* \in S_r$ given by $w^* =
(r+1 - w_{r}) (r+1 - w_{r-1}) \cdots (r+1 - w_1)$ has descent set
$D(w^*) = D(\alpha^*)$. If $u \in S_{r+l}$ is a shuffle of $w \in
S_r$ and $v \in S_l$, where $v \in S_l$ uses the letters
$r+1,r+2,\ldots, r+l$, then $u^*$ is a shuffle of $v^* \in S_l$ and
$w^* \in S_r$ where $w^* \in S_r$ uses the letters $l+1,l+2,\ldots,
r+l$.  Thus $\nu(L_{\C(v)}) \,\nu( L_{\C(w)}) =
L_{\C(v)^*}\,L_{\C(w)^*} = \nu(L_{\C(w)}\, L_{\C(v)})$, showing that
$\nu$ is an algebra map.  That $\nu$ is an involution is clear from
the definition.

The second statement can be deduced from the fact that $\nu$
commutes with the map $\alpha \mapsto \{\beta \mid D(\beta) \subset
D(\alpha)\}$.
\end{proof}

%Define (as in Example~\ref{ex:schur}) two functions
%$f^{\weak},f^{\strict}:\P \rightarrow \N \cup \{\infty\}$ by
%$f^{\weak}(n) = n$ and $f^{\strict}(n) = n-1$.
%
%\begin{proposition}
%Let $(P,O)$ be a finite $\T$-labelled poset.  Suppose $$O(s,t) \in
%\{f^{\weak},f^{\strict}\}$$ for each covering relation $s \lessdot
%t$.  Then $K_{P,O}(x)$ is a quasi-symmetric function.
%\end{proposition}
%
%A $\T$-labelled poset satisfying the conditions of the proposition
%is called {\it oriented} in~\cite{McN}.  Stanley's
%$(P,\omega)$-partitions are special cases of $(P,O)$-tableaux, for
%such posets.  If $f \in \QSym$ then $f$ is $m$-positive if and only
%if is a non-negative linear combination of the $M_\alpha$.

\section{Posets and $P$-partitions}
\label{sec:Ppart}

\subsection{Posets and cell transfer }
Let $(P, \leq)$ be a possibly infinite poset.  Let $s, t \in P$.
We say that $s$ {\it covers} $t$ and write $s \gtrdot t$ if for
any $r \in P$ such that $s \geq r \geq t$ we have $r = s$ or $r =
t$. The {\it Hasse diagram} of a poset $P$ is the graph with
vertex set equal to the elements of $P$ and edge set equal to the
set of covering relations in $P$.  If $Q \subset P$ is a subset of
the elements of $P$ then $Q$ has a natural induced subposet
structure. If $s, t \in Q$ then $s \leq t$ in $Q$ if and only if
$s \leq t$ in $P$. Call a subset $Q \subset P$ {\it connected} if
the elements in $Q$ induce a connected subgraph in the Hasse
diagram of $P$.

If $P$ and $Q$ are posets then the {\it disjoint sum} $P \oplus Q$
is the poset with the union $P \sqcup Q$ of elements, such that $a
\leq b$ in $P \oplus Q$ if either $a \leq b \in P$ or $a \leq b \in
Q$.

%The {\it dual} $P^*$ of a poset $P$ is the poset with the same
%elements as $P$ and partial order given by $p \leq p'$ in $P^*$ if
%and only if $p' \leq p$ in $P$.    The {\it product} $P \times Q$
%has elements given by $\{(p,q) \mid p \in P \,\text{and}\, q \in
%Q\}$ and inequalities $(p,q) \leq (p',q')$ if $p \leq p'$ in $P$ and
%$q \leq q'$ in $Q$.

An {\it order ideal} $I$ of $P$ is an induced subposet of $P$ such
that if $s \in I$ and $s \geq t \in P$ then $t \in I$.  A subposet
$Q \subset P$ is called {\it convex} if for any $s, t \in Q$ and $r
\in P$ satisfying $s \leq r \leq t$ we have $r \in Q$.
Alternatively, a convex subposet is one which is closed under taking
intervals.  A convex subset $Q$ is determined by specifying two
order ideals $J$ and $I$ so that $J \subset I$ and $Q = \{s \in I
\mid s \notin J\}$. We write $Q = I/J$.  If $s \notin Q$ then we
write $s < Q$ if $s < t$ for some $t \in Q$ and similarly for $s >
Q$. If $s \in Q$ or $s$ is incomparable with all elements in $Q$ we
write $s \sim Q$. Thus for any $s \in P$, exactly one of $s <Q$, $s
> Q$ and $s \sim Q$ is true.

\medskip

Let $Q$ and $R$ be two finite convex subposets of $P$.  Define the
{\it cell transfer operations} $\wedge$ and $\vee$ on the ordered
pair $(Q,R)$ by
\begin{equation}\label{eq:wedge}
Q \wedge R = \{s \in R \mid s < Q \} \cup \{s \in Q \mid s \sim R
\;\text{or}\; s < R\}
\end{equation}
and
\begin{equation}\label{eq:vee}
Q \vee R = \{s \in Q \mid s > R \} \cup \{s \in R \mid s \sim Q
\;\text{or}\; s > Q\}.
\end{equation}

\begin{lemma}[\cite{LP} Lemma 3.1]
The subposets $Q \wedge R$ and $Q \vee R$ are both convex subposets
of $P$. We have $(Q \wedge R) \cup (Q \vee R) = Q \cup R$ and $(Q
\wedge R) \cap (Q \vee R) = Q \cap R$.
\end{lemma}

The operations $\vee$, $\wedge$ are not commutative, and $Q \cap R$
is a convex subposet of both $Q \vee R$ and $Q \wedge R$.

%\remind{Complete this example.}
\begin{example}
\label{ex:young}
The poset $(\N^2,\leq)$ of (positive) points in a quadrant has cover
relations $(i,j) \gtrdot(i-1,j)$ and $(i,j) \gtrdot (i,j-1)$.  To
agree with the ``English'' notation for Young diagrams the first
coordinate $i$ increases as we go down while the second coordinate
$j$ increases as we go to the right.
The order ideals of $(\N^2 \leq)$ can be identified with Young
diagrams or alternatively with partitions.  
Let $\lambda = (4,1,1)$ and $\mu = (3,2)$ be two partitions interpreted as order ideals of $(\N^2,\leq)$.
Then applying the definitions (\ref{eq:wedge}) and (\ref{eq:vee}) 
above one can check that $\lambda
\wedge \mu = (3,1)$ and $\lambda \vee \mu = (4,2,1)$.
Figure~\ref{fig:ct} illustrates this example.

\begin{figure}
\begin{center}
\input{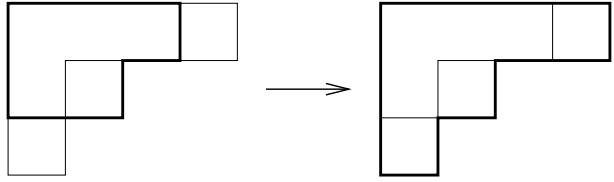}
\end{center}
 \label{fig:ct}
 \caption{Cell transfer for the Young shape posets $\lambda = (4,1,1)$ and $\mu = (3,2)$.}
  \end{figure}

 \end{example}

%Let $\P$ denote the set of positive integers and $\Z$ denote the
%set of integers. Let $\weak$ and $\strict$ be the two labels. A {\it labelling} $O$ of a poset $P$ is a map $O : \{(s,t)\in
%P^2 \mid s \gtrdot t \} \rightarrow \{ \weak, \strict \}$ labelling each edge
%$(s,t)$ of the Hasse diagram by one of the two labels.  A {\it $\T$-labelled
%poset} is an an ordered pair $(P,O)$ where $P$ is a poset, and $O$
%is a labeling of $P$. We shall refer to a labelled poset $(P,O)$ as $P$ when no
%ambiguity arises. The labeled posets are a particular case of $\T$-labelled posets introduced in \cite{LP}.

\subsection{Labeled posets}
Let $P$ be a poset.  A {\it labeling} $\theta$ of $P$ is an
injection $\theta: P \rightarrow \mathbb P$ into the positive
integers. A descent of the labeling $\theta$ of $P$ is a pair $p
\lessdot p'$ in $P$ such that $\theta(p) > \theta(p')$.  Let us say
that two labeled posets $(P,\theta_P)$, $(Q,\theta_Q)$ are
isomorphic if there is an isomorphism of posets $f: P \to Q$ so that
descents are preserved. That is if $p \lessdot p'$ then $\theta(p) <
\theta(p')$ if and only if $\theta'(f(p)) < \theta'(f(p'))$.  We say
that two labelings $\theta_1$ and $\theta_2$ of $P$ are {\it
equivalent} if the identity map on $P$ is an isomorphism of
$(P,\theta_1)$ and $(P,\theta_2)$.

\medskip
Let $(P,\theta)$ be a labeled poset.  If $Q \subset P$ is a
subposet, then it inherits a labeling $\theta|_Q$ by restriction.
When no confusion can arise, we will often denote $\theta|_Q$ simply
by $\theta$. Note however, that the descents of $\theta|_Q$ are not
completely determined by the descents of $\theta$, unless $Q$ is a
convex subset of $P$.

%If $\theta$ is a labeling of a finite poset $P$, then a conjugate
%labeling $\theta'$ is given by $\theta'(p) = m - \theta(p)$ where
%$m$ is a sufficiently large integer.  We then say $(P,\theta)$ is
%conjugate to $(P,\theta')$.

\medskip

Let $(P,\theta_P)$ and $(Q,\theta_Q)$ be labeled posets.  Then the
{\it disjoint sum} $(P \oplus Q, \theta^\oplus)$ is the labeled
poset (defined up to equivalence of labelings) where $\theta^\oplus$
has the same descents as the function
\[
f(a) = \begin{cases} \theta_P(a) & \mbox{if $a \in P$,
and} \\
\theta_Q(a) & \mbox{if $a \in Q$.} \end{cases}
\]
%The {\it dual} of a labeled poset $(P,\theta)$ is the labeled poset
%$(P^*,\theta^*)$ where $P^*$ is the dual of $P$ and $\theta =
%\theta^*$ as functions.
%The {\it product} $(P \times Q, \theta^\times)$ is the labeled poset
%(defined up to equivalence of labelings) where $\theta^\times$ has a
%descent at $(p,q) \lessdot (p',q')$ if and only if either $q = q'$
%and $p \lessdot p'$ is a descent of $\theta_P$ or $p = p'$ and $q
%\lessdot q'$ is a descent of $\theta_Q$.

%\remind{Complete this example}
\begin{example}
Let $P$ be the diamond poset with elements $a < b,c < d$ and labeling $\theta_P$ given by $\theta_P(a) = 2$, $\theta_P(b) = 1$, $\theta_P(c) = 3$, and $\theta_P(d) = 4$.  Let $Q$ be the chain with elements $e < f < g$ and labeling $\theta_Q$ given by $\theta_Q(e) = 1$, $\theta_Q(f) = 3$, and $\theta_Q(g) = 2$.  The one possible labeling $\theta^\oplus$ for the disjoin sum $P \oplus Q$ is given by 
$\theta^\oplus(a,b,c,d,e,f,g) = 4,3,5,7,1,6,2$.
In Figure \ref{lct1}, the three labeled posets $(P,\theta_P)$, $(Q,\theta_Q)$, and $(P \oplus Q, \theta^\oplus)$ are shown.  Note that we have some freedom in choosing the labelling $\theta^\oplus$.

\begin{figure}[!ht]

\begin{center}
\input{lct1.pstex_t}
\end{center}\caption{}
\label{lct1}
\end{figure}
 
\end{example}

\subsection{$P$-partitions}
\begin{definition}
A {\it $(P, \theta)$-partition} is a map $\sigma: P \to \mathbb P$
such that for each covering relation $s \lessdot t$ in $P$ we have
\begin{align*}
\sigma(s) &\leq \sigma(t) & \mbox{ if $\theta(s) < \theta (t)$,} \\
\sigma(s) &< \sigma(t) & \mbox{ if $\theta(t) < \theta (s)$.}
\end{align*}
If $\sigma: P \to \mathbb P$ is any map, then we say that $\sigma$
{\it respects} $\theta$ if $\sigma$ is a $(P, \theta)$-partition.
\end{definition}

Denote by $\mathcal{A}(P,\theta)$ the set of all
$(P,\theta)$-partitions.  Clearly $\mathcal{A}(P,\theta)$ depends on
$(P,\theta)$ only up to isomorphism.  If $P$ is finite then one can
define the formal power series $K_{P,\theta}(x_1,x_2,\ldots) \in
\Z[[x_1,x_2,\ldots]]$ by
\[
K_{P,\theta}(x_1,x_2,\ldots) = \sum_{\sigma \in
\mathcal{A}(P,\theta)} x_1^{\# \sigma^{-1}(1)} x_2^{\#
\sigma^{-1}(2)} \cdots.
\]
The composition $\wt(\sigma) = (\# \sigma^{-1}(1), \#
\sigma^{-1}(2), \ldots)$ is called the {\it weight}§ of $\sigma$.

Let $P$ be a poset with $n$
elements.  Recall that a linear extension of $P$ is a bijection $e:
P \to \{1,2,\ldots,n\}$ satisfying $e(x) \leq e(y)$ if $x \leq y$ in
$P$.  The Jordan-Holder set $\J(P,\theta)$ of $(P,\theta)$ is the
set $$\{\theta(e^{-1}(1)) \theta(e^{-1}(1)) \cdots
\theta(e^{-1}(n))\mid \text{$e$ is a linear extension of $P$}\}.$$ It
is a subset of the set $\SG(\theta(P))$ of permutations of
$\theta(P) \subset \mathbb P$.

\begin{example}
Suppose $C$ is a chain $c_1 < c_2 < \ldots < c_n$ with $n$ elements
and $w = w_1w_2\dots w_n \in \SG_n$ a permutation of
$\{1,2,\ldots,n\}$. Then $(C,w)$ can be considered a labeled poset,
where $w(c_i) = w_i$.  In this case we have $K_{C,w} = L_{\C(w)}$.
\end{example}

\begin{theorem}[\cite{Sta}]
\label{thm:Ppart} The generating function $K_{P,\theta}$ is
quasi-symmetric. We have $K_{P,\theta} = \sum_{w \in \J(P,\theta)}
L_{D(w)}$.
\end{theorem}

In particular, $K_{P,\theta}$ is $L$-positive.  This motivates our
treatment of $K_{P,\theta}$ as ``skew''-analogues of the functions
$L_\alpha$. Let $Q$ and $R$ be two finite convex subposets of
$(P,\theta)$.
\begin{theorem}[\cite{LP}]
\label{thm:celltransfer} The difference $K_{Q \wedge
R,\theta}K_{Q\vee R,\theta} - K_{Q,\theta}K_{R,\theta}$ is
$M$-positive.
\end{theorem}
The main theorem of~\cite{LP} generalizes
Theorem~\ref{thm:celltransfer} to more general labelings.  We will return to a discussion 
of these more general labelings in Section~\ref{sec:gen}.

\begin{example}
Let $P = \lambda$ be the poset of squares in the Young diagram of a
partition $\lambda$ as in Example~\ref{ex:young}.  Let $\theta_{\rm reading}$ 
be the labeling of $\lambda$ obtained from the bottom to top
row-reading order.  Then $K_{\lambda,\theta_{\rm reading}}$ is equal
to the Schur function $s_\lambda$. In~\cite{LP} it is conjectured
and in~\cite{LPP} it is shown that in this case the expression of
Theorem~\ref{thm:celltransfer} is Schur positive, which implies
monomial positivity.
\end{example}

\section{Cell transfer for $P$-partitions}
By Theorem~\ref{thm:Ppart}, the expression $K_{Q \wedge
R,\theta}K_{Q\vee R,\theta} - K_{Q,\theta}K_{R,\theta}$ is always a
quasi-symmetric function.  We now show that this difference is
$L$-positive.

Let $(P,\theta)$ be a labeled poset and let $Q$ and $R$ be convex
subsets.  In~\cite{LP}, we gave a weight preserving injection
$$
\eta: \A(Q,\theta) \times \A(R,\theta) \longrightarrow \A(Q \wedge
R,\theta) \times \A(Q\vee R, \theta).
$$
The injection $\eta$ satisfies additional crucial properties.  First
let us say that $i \neq  j$ are {\it adjacent} in a multiset $T$ (of
integers) if $i,j \in T$ and for any other $t \in T$ both $i \leq t
\leq j$ and $ j \leq t \leq i$ fail to hold.

\begin{proposition}
\label{prop:ctprop} Suppose $\omega \in \A(Q,\theta)$ and $\sigma
\in \A(R,\theta)$ and $\eta(\omega,\sigma) = (\omega \wedge \sigma,
\omega \vee \sigma)$. Let $p \in Q \cup R$.
\begin{enumerate}
\item
If $p \in Q \cap R$, then $\{\omega(p),\sigma(p)\} =
\{\omega\wedge\sigma(p),\omega\vee\sigma(p)\}$.  Furthermore,
suppose $\omega(p)$ and $\sigma(p)$ are adjacent in the multiset
$\omega(Q) \cup \sigma(R)$.  Then $\omega \wedge \sigma (p) =
\omega(p)$ and $\omega \vee \sigma(p) = \sigma(p)$.
\item
If $p \in Q \wedge R$ but $p \notin Q \cap R$ then $\omega \wedge
\sigma(p) = \omega(p)$ if $p \in Q$ and $\omega \wedge \sigma(p) =
\sigma(p)$ if $p \in R$.
\item
If $p \in Q \vee R$ but $p \notin Q \cap R$ then $\omega \vee
\sigma(p) = \omega(p)$ if $p \in Q$ and $\omega \vee \sigma(p) =
\sigma(p)$ if $p \in R$.
\end{enumerate}
\end{proposition}
Roughly speaking, Proposition~\ref{prop:ctprop}(1) says that if $p
\in Q \cap R$, then one obtains $(\omega \wedge \sigma(p),\omega
\vee \sigma(p))$ by possibly ``swapping'' $\omega(p)$ with
$\sigma(p)$; in addition, no swapping occurs if $\omega(p)$ and
$\sigma(p)$ are adjacent in $\omega(Q) \cup \sigma(R)$.

\begin{proof}
Let $S \subset Q \cap R$.  In~\cite{LP}, $(\omega \wedge \sigma)_S:
Q \wedge R \rightarrow \P$ was defined by $$ (\omega \wedge
\sigma)_S(x) =
\begin{cases}
% \omega(x) & \mbox{if $x
%\in Q$ but $x \notin R$,} \\
\sigma(x) & \mbox{if $x \in R \backslash Q$ or $x \in S$,} \\
\omega(x) & \mbox{otherwise,}
\end{cases}
$$
and $(\omega \vee \sigma)_S: Q \vee R \rightarrow \P$ by
$$
(\omega \vee \sigma)_S(x) = \begin{cases}
%\sigma(x) & \mbox{if $x
%\in R$ but $x \notin Q$,} \\
\omega(x) & \mbox{if $x \in Q \backslash R$ or $x \in S$,} \\
\sigma(x) & \mbox{otherwise.}
\end{cases}
$$
All statements except the last sentence of (1) follows from the
definition of $\eta(\omega,\sigma)$ as $((\omega \wedge
\sigma)_{S^\diamond},(\omega \vee \sigma)_{S^\diamond})$ for some
choice of $S = S^\diamond$ defined in the proof of
Theorem~\ref{thm:celltransfer} in~\cite{LP}.  The last statement of
(1) follows from the fact that $S^\diamond$ is defined to be the
{\it smallest} set such that $((\omega \wedge
\sigma)_{S^\diamond},(\omega \vee \sigma)_{S^\diamond})$ is an
element of $\A(Q \wedge R,\omega) \times \A(Q\vee R, \omega)$.  More
precisely, if $p \in Q \cap R$ is such that $\omega(p)$ and
$\sigma(p)$ are adjacent then $p \notin S^\diamond$.
\end{proof}

Now consider the labeled posets $(Q \oplus R,\theta^\oplus)$ and
$((Q \vee R) \oplus (Q \wedge R), \theta^{\vee\wedge})$, where we
shall pick $\theta^\oplus$ and $\theta^{\vee\wedge}$ as follows.

For each $p \in Q \cap R$, we ``duplicate'' $\theta(p)$ by picking
$\theta(p)' > \theta(p)$ so that for every $x \in Q \cup R$ such
that $x \neq p$ we have $\theta(p)' < \theta(x)$ if and only if
$\theta(p) < \theta(x)$; also the duplicates satisfying the same
inequalities as the originals so that $\theta(p)' < \theta(x)'$ if
and only if $\theta(p) < \theta(x)$. This describes a total order on
$\{\theta(p) \mid p \in Q \cup R\} \cup \{\theta(p)' \mid p \in Q
\cap R\}$.  Note that we may need to replace $\theta$ with an
equivalent labeling so that there is enough ``space'' to insert the
primed letters.

Now suppose $p \in Q \cap R$.  Denote the copy of $p$ inside $Q
\subset Q \oplus R$ by $p_Q$ and the copy of $p$ inside $R \subset Q
\oplus R$ by $p_R$.  Similarly, denote the elements of $(Q \vee R)
\oplus (Q \wedge R)$.  We define
$$
\theta^\oplus(p) = \begin{cases} \theta(p) & \mbox{if $p
\notin Q \cap R$,} \\
\theta(p)' & \mbox{if $p = p_Q$,} \\
\theta(p) & \mbox{if $p = p_R$} \end{cases}
$$
and
$$
\theta^{\vee\wedge}(p) = \begin{cases} \theta(p) & \mbox{if $p
\notin Q \cap R$,} \\
\theta(p)' & \mbox{if $p = p_{Q\wedge R}$,} \\
\theta(p) & \mbox{if $p = p_{Q \vee R}$.} \end{cases}
$$
Clearly the descents of $\theta^\oplus$ (or $\theta^{\vee\wedge}$)
on either component agree with the descents of that component as a
convex subposet of $(P,\theta)$.

\begin{theorem}
\label{thm:main} The difference $K_{Q \wedge R,\theta}K_{Q\vee
R,\theta} - K_{Q,\theta}K_{R,\theta}$ is $L$-positive.
\end{theorem}
\begin{proof}
Let $|Q| + |R| = n = |Q\vee R| + |Q \wedge R|$ and suppose $\alpha:
Q \oplus R \rightarrow [n]$ is a linear extension.  Then $\alpha$ in
particular gives an element $(\alpha|_Q,\alpha|_R)$ of $\A(Q,\theta)
\times \A(R,\theta)$. Using Proposition~\ref{prop:ctprop}, we see
that $\eta(\alpha|_Q,\alpha|_R) = (\beta|_{Q \wedge R}, \beta|_{Q
\vee R})$ arises from a linear extension $\beta: (Q \wedge R) \oplus
(Q \vee R) \rightarrow [n]$ (in other words the union $\beta|_{Q
\wedge R}\cup \beta|_{Q \vee R}$ is exactly the interval $[n]$).

We claim that the two words
\begin{align*}
a_\alpha & =
\theta^\oplus(\alpha^{-1}(1))\theta^\oplus(\alpha^{-1}(2))\dots
\theta^\oplus(\alpha^{-1}(n)) \\ b_\beta &=
\theta^{\vee\wedge}(\beta^{-1}(1))\theta^{\vee\wedge}(\beta^{-1}(2))\dots
\theta^{\vee\wedge}(\beta^{-1}(n))
\end{align*}
have the same descent set.  Again by Proposition~\ref{prop:ctprop},
the word $b = b_1 b_2 \dots b_n$ is obtained from $a = a_1a_2\dots
a_n$ by swapping certain pairs $(a_i,a_j)$ where $a_i =
\theta^\oplus(p_Q)$ and $a_j = \theta^\oplus(p_R)$ for some $p \in Q
\cap R$.

By definition $\theta^\oplus(p_Q) = \theta^{\vee\wedge}(p_{Q \wedge
R})$ and $\theta^\oplus(p_R) = \theta^{\vee\wedge}(p_{Q \vee R})$ so
swapping occurs if and only if $(\alpha(p_Q),\alpha(p_R)) =
(\beta(p_{Q \vee R}),\beta(p_{Q \wedge R}))$.  By the last statement
of Proposition~\ref{prop:ctprop}~(1), this never happens if
$\alpha|_Q(p_Q)$ and $\alpha|_R(p_R)$ are adjacent in $[n]$, which
is equivalent to $|i-j| = 1$.  Thus swapping $(a_i,a_j)$ is the same
as swapping a pair of non-neighboring letters
$(\theta(p),\theta(p)')$ in the word $a_1a_2 \dots a_n$, which
preserves descents by our choice of $\theta(p)'$.

We have $K_{Q,\theta}K_{R,\theta} = \sum_{\alpha} L_{D(a_\alpha)}$
and $K_{Q \wedge R,\theta}K_{Q \vee R,\theta} = \sum_{\beta}
L_{D(b_\beta)}$, where the summations are over linear extensions of
$Q \oplus R$ and $(Q \wedge R) \oplus (Q \vee R)$.  Since $\eta$
induces an injection from the first set of linear extensions to the
second, we conclude that $K_{Q \wedge R,\theta}K_{Q\vee R,\theta} -
K_{Q,\theta}K_{R,\theta}$ is $L$-positive.
\end{proof}

\begin{example}
Let $P$ be the poset on the $5$ elements $A, B, C, D, E$ given by  the cover
relations $A<B$, $A<C$, $B<D$, $B<E$, $C<D$, $C<E$. Take the
following labeling $\theta$ of $P$: $\theta(A)=2$, $\theta(B)=1$,
$\theta(C)=4$, $\theta(D)=5$, $\theta(E)=3$. Take the two ideals $Q =
\{A,B,C,D\}$, $R= \{A,B,C,E\}$ of $P$. Form the disjoint sum poset
$Q \oplus R$.  The elements $A, B,C \in Q \cap R$ have two images in the
newly formed poset: $A_Q,B_Q,C_Q$ and $A_R,B_R,C_R$. The labels of
$Q \oplus R$ are formed according to the rule above: for $X = A,B,C$
we have $\theta^{\oplus}(X_Q)=\theta(X)'$ while
$\theta^{\oplus}(X_R)=\theta(X)$. The resulting labeling is shown in
Figure \ref{lct2}, with $\theta^{\oplus}$ taking the values $\{1 < 1' < 2 < 2' < 3 < 4 < 4' <5\}$ .

\begin{figure}[!ht]
\begin{center}
\input{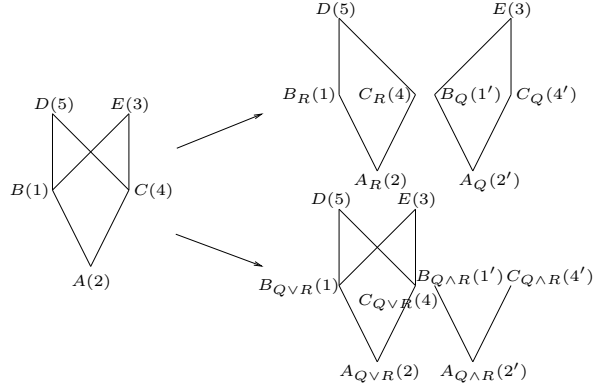}
\end{center}
\caption{The labelings $\theta^{\oplus}$ of $Q \oplus R$ and $\theta^{\wedge \vee}$ of $(Q \wedge R) \oplus (Q \vee R)$ formed from a labeling $\theta$ of $P$.  Labels are shown in parentheses.}
\label{lct2}
\end{figure}

Similarly, we obtain the labeling $\theta^{\wedge \vee}$ of $(Q \wedge
R) \oplus (Q \vee R)$, as shown on Figure \ref{lct2}. Clearly each edge in the
Hasse diagrams of $Q \oplus R$ and $(Q \wedge R) \oplus (Q \vee R)$
is a descent if and only if it is in the Hasse diagram of $P$.

Now, to illustrate the proof of Theorem~\ref{thm:main} take a particular extension of
$Q \oplus R$, namely $\alpha$ defined by $\alpha^{-1}([8])=(A_Q,
A_R, B_Q, C_Q,E,B_R,C_R,D)$. Performing cell transfer we get
$\beta = \eta(\alpha)$ with $S^\diamond = \{B,C\}$ in the notation of the proof of 
Proposition~\ref{prop:ctprop}, so that
\[\beta^{-1}([8]) =
(A_{Q \wedge R},A_{Q \vee R},B_{Q \vee R},C_{Q \vee R},E,B_{Q \wedge R},C_{Q \wedge R},D)\] 
In this case
$a_{\alpha}=(2',2,1',4',3,1,4,5)$ and $b_{\beta}=(2',2,1,4,3,1',4',5)$.
The pairs that got swapped are $(1,1')$ and $(4,4')$.  Note also that the pair $(2,2')$
did not get swapped, which we know cannot happen since those labels
are neighbors are in the word $a_\alpha$.  It is clear that the descents
in $a_{\alpha}$ are indeed the same as in $b_{\beta}$.

\begin{figure}[h!]
\label{lct3}
\begin{center}
\input{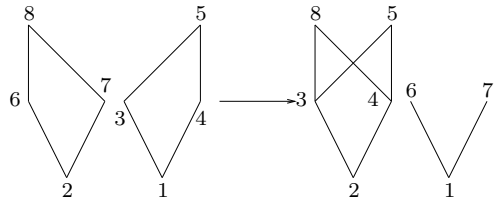}
\end{center}
\caption{The linear extension $\alpha$ of $Q \oplus R$ and the linear extension $\beta$ of
$(Q \wedge R) \oplus (Q \vee R)$ obtained by cell transfer.}
\end{figure}

\end{example}

Comparing Theorem~\ref{thm:main} and Theorem~\ref{thm:Ppart}, we
obtain the following question.
\begin{question}
\label{q:main} When is the difference $K_{Q \wedge R,\theta}K_{Q\vee
R,\theta} - K_{Q,\theta}K_{R,\theta}$ itself of the form $K_{S,\xi}$
for some labeled poset $(S,\xi)$?
\end{question}
In other words, we are asking for another (hopefully natural)
operation $\sharp$ on convex subsets $Q$ and $R$ of a labeled poset
$(P,\theta)$ so that
$$
K_{Q \sharp R, \theta^\sharp} = K_{Q \wedge R,\theta}K_{Q\vee
R,\theta} - K_{Q,\theta}K_{R,\theta}.
$$
We will give an affirmative answer to Question~\ref{q:main} for the case of chains in
Section~\ref{sec:wave}.  As the following example shows, 
the answer to Question~\ref{q:main} is not affirmative in general.

\begin{example}
Let $P$ be the poset with four elements $\{a,b,c,d\}$ and relations $a < b, a < c, a< d$.
Give $P$ the labeling $\theta(a) = 4$, $\theta(b) = 1$, $\theta(c) = 2$, and $\theta(d) = 3$. 
Let $Q$ be the ideal $\{a,b\}$ and $R$ be the ideal $\{a,c,d\}$.  Then the difference 
$K_{Q \wedge R,\theta}K_{Q\vee
R,\theta} - K_{Q,\theta}K_{R,\theta}$ is given by 
\begin{equation}
\label{eq:bad}
d = L_{1}(L_{1111} + 2L_{112} + 2L_{121} + L_{13}) - L_{11} (L_{12} + L_{111}).
\end{equation}
We will argue that $d$ is not equal to $K_{S,\theta_S}$ for any $(S,\theta_S)$.
First we claim that no term $L_{\alpha}$ in the $L$-expansion of $d$ has $\alpha_1 > 1$.
It is not difficult to see directly from the shuffle product that the expansion of each term in $d$ 
has six $L_{\alpha}$ terms with $\alpha_1 > 1$ (in fact $\alpha_1 = 2$) and these cancel out
by Theorem~\ref{thm:main}.  

Thus using Theorem~\ref{thm:Ppart} we conclude that if $d = K_{S,\theta_S}$ then $S$ must be a 
five element poset with a unique minimal element.  Also one computes from (\ref{eq:bad})
that $S$ must have exactly 10 linear 
extensions.  No poset $S$ has these properties.
\end{example}

\begin{remark}
By carefully studying the cell transfer injection $\eta$ of~\cite{LP}, 
one can also give an affirmative answer to Question~\ref{q:main} for the case
where $P$ is a tree, and $Q$ and $R$ are order ideals so that both $Q/Q\cap R$ and
$R/Q \cap R$ are connected.
\end{remark}

\begin{remark}
Question~\ref{q:main} can be asked for the $\T$-labeled posets
of~\cite{LP} and also for the differences of products of skew Schur
functions studied in~\cite{LPP}.  However, we will not investigate
these questions in the current article.
\end{remark}

\section{Chains and fundamental quasi-symmetric functions} \label{sec:chains}
\subsection{Cell transfer for compositions}
Let $(C_n,w)$ be the labeled chain corresponding to the permutation
$w \in \SG_n$. Let us consider $C_n$ to consist of the elements
$\{c_1 < c_2 < \dots < c_n\}$, so that $w:C_n \to \P$ is given by
$w(c_i) = w_i$. The convex subsets $C[i,j]$ of $C_n$ are in
bijection with intervals $[i,j] \subset [n]$.

Let $Q = [a,b]$ and $R = [c,d]$ and assume that $a \leq c$. Then
we have the following two cases:
\begin{enumerate}
\item If $b \leq d$ then $Q \wedge R = Q$ and $Q \vee R = R$.
\item If $b \geq d$ then $Q \wedge R = [a,d]$ and $Q \vee R = [b,c]$.
\end{enumerate}
Thus to obtain a non-trivial cell transfer we assume that $a < c
\leq d < b$.  Let $w[i,j]$ denote the word $w_i w_{i+1} \dots
w_{j}$.  Theorem~\ref{thm:main} then says that the difference
\begin{equation}
\label{eq:Lchain} L_{\C(w[a,d])} \, L_{\C(w[c,b])} - L_{\C(w[a,b])}
\,L_{\C(w[c,d])}
\end{equation}
is $L$-positive.

\begin{figure}[h!]
\begin{center}
\input{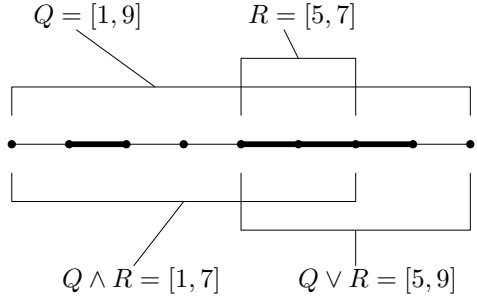}
\caption{An example of the cell transfer operation for chains, here $w=(2,1,5,4,3,7,8,9,6)$, $a=1, b=9, c=5, d=7$.} \label{fig:chains}
\end{center}
\end{figure}

We now make the difference (\ref{eq:Lchain}) more precise by
translating into the language of compositions and descent sets.  Let
$\alpha = (\alpha_1,\alpha_2,\ldots,\alpha_k)$ and $\beta =
(\beta_1,\beta_2,\ldots,\beta_l)$ be an ordered pair of
compositions.  Say $\beta$ can be {\it found inside} $\alpha$ if
there exists a non-negative integer $m \in [0,|\alpha|-|\beta|]$ so
that $D(\beta) + m$ coincides with $D(\alpha)$ restricted to
$[m+1,m+|\beta|-1]$.  We then say that $\beta$ can be found inside
$\alpha$ at $m$.  A composition can be found inside another in many
different ways. For example if $\beta = (1)$ then one may pick $m$
to be any integer in $[0,|\alpha|-1]$.

Now for a composition $\alpha = (\alpha_1,\alpha_2,\ldots,\alpha_k)
\vdash n$ and an integer $x \in [1,n]$ we define two new
compositions $\alpha^{x\leftarrow}, \alpha^{x\to} \vdash x$ as
follows. We define $\alpha^{x\leftarrow} =
(\alpha_1,\alpha_2,\ldots,\alpha_{r-1},a)$ where $a,r$ are the
unique integers satisfying $1 \leq a \leq \alpha_r$ and $\alpha_1 +
\alpha_2 + \dots + \alpha_{r-1} + a = x$. Similarly, define
$\alpha^{x\to} = (b,\alpha_{s+1},\ldots,\alpha_k)$ where $b,s$ are
the unique integers satisfying $1 \leq b \leq \alpha_s$ and $b +
\alpha_{s+1} + \dots + \alpha_k = x$.  If $\beta$ can be found
inside $\alpha$ at $m$, we set $\alpha \wedge_m \beta = \alpha^{(m +
|\beta|)\leftarrow}$ and $\alpha \vee_m \beta =
\alpha^{(|\alpha|-m)\to}$.

The $L$-positive expressions in (\ref{eq:Lchain}) give the following
theorem.

\begin{theorem}
\label{thm:lp} Let $\alpha$ and $\beta$ be compositions such that
$\beta$ can be found inside $\alpha$ at $m$.  Then the difference
$$
L_{\alpha \wedge_m \beta}\, L_{\alpha \vee_m \beta} - L_\alpha \,
L_\beta
$$
is $L$-positive.
\end{theorem}

\begin{example}
Let us take chain $(C_9,w)$ with $w=(2,1,5,4,3,7,8,9,6)$ and $Q = [1,9]$ and $R = [5,7]$ so that
$a=1$, $b=9$, $c=5$, $d=7$. Then we get the situation shown in Figure \ref{fig:chains}, the thinner edges indicating descents. If $\alpha = (1,2,1,4,1)$ and $\beta = (3)$ then there are two ways to find $\beta$ inside $\alpha$, and Figure \ref{fig:chains} shows the way to find it at $m=5$. In this case Theorem \ref{thm:lp} says that $L_{(1,2,1,3)} L_{(4,1)}-L_{(1,2,1,4,1)} L_{(3)}$ is $L$-positive.
\end{example}

\begin{remark}
The operation $(\alpha,\beta) \mapsto (\alpha \wedge_m \beta, \alpha
\vee_m \beta)$ interacts well with the involutions $\nu$ and
$\omega$ of $\QSym$.  More precisely, if $\beta$ can be found inside
$\alpha$ at $m$ then $\beta^*$ can be found inside $\alpha^*$ at $m$
and $\bar \beta$ can be found inside $\alpha^*$ at
$|\alpha|-|\beta|-m$.
\end{remark}

\subsection{The $L$-positivity poset}
Fix a positive integer $n$.  Now define a poset structure
$(PC_n,\leq)$ (``Pairs of Compositions'') on the set $PC_n$ of
unordered pairs $\{\alpha,\beta\}$ of compositions satisfying
$|\alpha|+|\beta| = n$ by letting $\{\alpha,\beta\} \leq
\{\gamma,\delta\}$ if $L_\gamma \,L_\delta - L_\alpha \,L_\beta$ is
$L$-nonnegative.  The following result relies on factorization properties of $\QSym$
which we prove in Section~\ref{sec:fact}.

\begin{proposition}
The relation $\{\alpha,\beta\} \leq \{\gamma,\delta\}$ if $L_\gamma
\,L_\delta \geq_L L_\alpha \,L_\beta$ defines a partial order on the
set $PC_n$.
\end{proposition}
\begin{proof}
Reflexivity and transitivity of $\leq$ are clear. Suppose we have
both $\{\alpha,\beta\} \leq \{\gamma,\delta\}$ and $
\{\gamma,\delta\} \leq \{\alpha,\beta\}$ then we must have $L_\gamma
\,L_\delta = L_\alpha \,L_\beta$.  By Corollary~\ref{cor:UFD} and
Proposition~\ref{prop:irred} we must have $\{\alpha,\beta\} =
\{\gamma,\delta\}$.  Thus $\leq$ satisfies the symmetry condition of
a partial order.
\end{proof}

For an unordered pair of compositions $\{\alpha,\beta\}$ we
unambiguously define another unordered pair $\{\alpha \vee \beta,
\alpha \wedge \beta\}$ as follows.  Suppose $|\alpha| \geq |\beta|$.
If $\beta$ can be found inside of $\alpha$, we pick the smallest $m
\in (0, |\alpha|-|\beta|)$ where this is possible and set $\alpha
\wedge \beta = \alpha \wedge_m \beta$ and $\alpha \vee \beta =
\alpha \vee_m \beta$. Otherwise we set $\{\alpha \vee \beta, \alpha
\wedge \beta\} = \{\alpha,\beta\}$.

\begin{conjecture}
\label{conj:max} The maximal elements of $PC_n$ are exactly the
pairs $\{\alpha,\beta\}$ for which $\{\alpha,\beta\} =
\{\alpha\wedge\beta,\alpha \vee \beta\}$.
\end{conjecture}

Note that Conjecture~\ref{conj:max} is compatible with the two
involutions $\omega$ and $\nu$ of $\QSym$.

\begin{remark}
(i) Conjecture~\ref{conj:max} has been verified by computer up to
$n=10$.

(ii) A result similar to Conjecture~\ref{conj:max} holds for the
case of Schur functions: the pairs of partitions corresponding to
Schur-maximal products $s_\lambda \, s_\mu$ are exactly those
partitions fixed by ``skew cell transfer''; see~\cite{LPP2}.
\end{remark}

\begin{example}
In Figure \ref{lct5} the poset $PC_4$ is shown, a composition $\alpha$ being represented by a 
chain $(C,w)$ satisfying $K_{C,w} = L_\alpha$.  The actual labeling $w$ of the chain is not shown,
instead the descents of $w$ are marked with thin edges. The elements of the bottom row are single compositions of size $4$ since the second composition in this case is empty. 

\begin{figure}[h!]
\begin{center}
\input{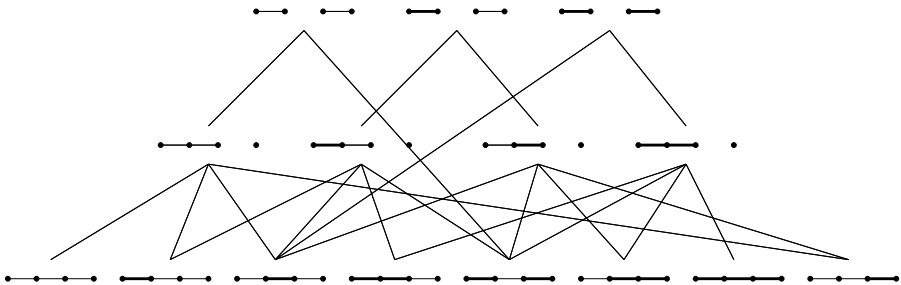}
\caption{Partial order $PC_4$ on pairs of compositions, descents are drawn as thin edges.} \label{lct5}
\end{center}
\end{figure}

One can see that maximal elements are exactly the ones for which one of the two compositions cannot be found inside the other. In this case those are exactly pairs $(\alpha, \beta)$ such that $|\alpha|=|\beta|=2$.

\end{example}

\section{Wave Schur functions}
\label{sec:wave}
In this section we define new generating functions called wave Schur functions.  
We first show that they are $L$-positive, and then prove a determinantal formula for them.

\subsection{Wave Schur functions as $P$-partition generating functions}

%Define a poset $QPos$ on pairs of compositions with the same total
%degree given by $L$-positivity of the products of $L$-quasi
%symmetric functions.

%\begin{conjecture}
%The maximal elements of $QPos$ are exactly the pairs such that cell
%transfer cannot be performed.  The number of such elements is given
%by ...
%\end{conjecture}

%It is not the case that this poset $QPos$ is generated by the
%inequalities given by cell transfer (letting one $L_\alpha$ have
%degree 1).  It is easy to see that two products with the same
%multi-degree (i.e. if ${\rm deg}(L_\alpha) = {\rm deg}(L_\delta)$
%and ${\rm deg}(L_\beta) = {\rm deg}(L_\gamma)$ are incomparable.

\begin{figure}
\begin{center}
\input{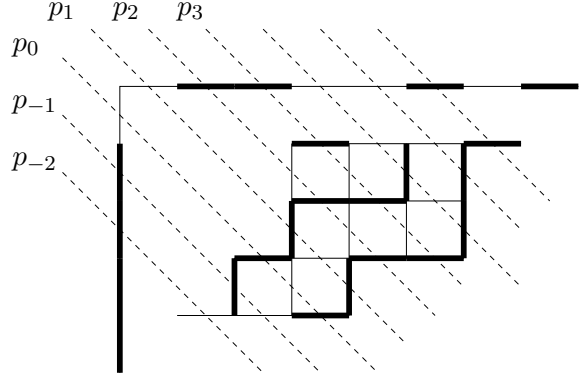}
\end{center}
\caption{An edge labeling $O_\p$ of $(\N^2,\leq)$.} \label{fig:gens}
\end{figure}

The poset $(\N^2,\leq)$ of (positive) points in a quadrant has cover
relations $(i,j) \gtrdot(i-1,j)$ and $(i,j) \gtrdot (i,j-1)$.  To
agree with the ``English'' notation for Young diagrams the first
coordinate $i$ increases as we go down while the second coordinate
$j$ increases as we go to the right.  Let us fix a sequence of
``strict--weak'' assignments $\p = \{p_i \in \{{\weak},{\strict}\}
\mid i \in \Z\}$. Let $\overline{\weak} = \strict$ and
$\overline{\strict} = \weak$. Define an {\it edge-labeling} (or {\it
orientation} in the language of~\cite{McN}) $O_\p$ as a function
from the covers of $\N^2$ to $\{\weak,\strict\}$ by
\begin{align*}
O_\p((i,j) \gtrdot (i-1,j)) &= \overline{p_{j-i+1}} & \text{and} &&
O_\p((i,j) \gtrdot (i,j-1)) &= p_{j-i}.
\end{align*}
An example of an such an edge-labeling $O_\p$ is given in
Figure~\ref{fig:gens}, where
$$\ldots p_{-3}, p_{-2}, p_{-1}, p_{0}, p_{1}, p_{2}, p_{3},
p_{4}, p_{5}, p_{6}, p_{7}, p_{8}, \ldots =$$ $$\ldots \strict,
\strict, \strict, \weak, \strict, \weak, \weak, \strict, \strict,
\weak,\strict,\weak \ldots.$$ The lines show the diagonals along
which $O_\p$ alternates between weak and strict edges.  We have
labeled weak edges thick and strong edges thin (agreeing with the way we
labeled chains in Section~\ref{sec:chains}).

%\remind{label the figure with $p_i$ to show where $p_0$ is}
%is a descent of $\theta_\p$ if and only if $p_{j-i} = \strict$ and
%$(i,j) \gtrdot (i,j-1)$ is a descent if and only if
%$\overline{p_{j-i-1}} = \strict$

In the following definition, $\lambda/\mu$ denotes a Young diagram
$\{(i,j) \mid \mu_i \leq j \leq \lambda_i\}$ considered as a
subposet of $(\N^2,\leq)$.
\begin{definition}
A {\it wave $\p$-tableau} of shape $\lambda/\mu$ is a function $T:
\lambda \to \P$ such that for each cover $s \lessdot t$ we have
\begin{align*}
T(s) &< T(t) & \mbox{if $O_\p(s \lessdot t) = \strict$,} \\
T(s) &\leq T(t) & \mbox{if $O_\p(s \lessdot t) = \weak$.}
\end{align*}
The wave Schur function $s^\p_{\lambda/\mu}$ is given by the weight
generating function
$$s^\p_{\lambda/\mu}(x_1,x_2,\ldots) = \sum_T x_1^{\#T^{-1}(1)}  x_2^{\#T^{-1}(2)} \dots$$
of all wave $\p$-tableaux of shape $\lambda/\mu$.
\end{definition}

The {\it standard} ``strict-weak'' assignment is given by $\p =
\{p_i\}$ where $p_i = \weak$ for all $i$.  In this case a wave
$\p$-tableau is a usual semistandard tableau, and the wave Schur
function is the usual Schur function.  Note, however, that in
general a wave Schur function is not symmetric.  However, wave Schur functions are always $(P,\theta)$-partition generating
functions.
\begin{proposition} \label{glp}
Let $\lambda/\mu$ be a skew shape.  There exists a (vertex) labeling
$\theta_\p: \lambda/\mu \to \P$ such that $(s \lessdot t)$ is a
descent of $\theta_\p$ if and only if $O_\p(s \lessdot t) =
\strict$.  Thus $s_{\lambda/\mu}^\p = K_{\lambda/\mu,\theta_\p}$.
\end{proposition}

\begin{proof}
We shall prove the result by induction on the number of boxes in
$\lambda/\mu$.  Let $(i,j)$ be any outer corner of $\lambda/\mu$. In
other words there are no boxes to the bottom right of $(i,j)$, and
if we remove $(i,j)$ from $\lambda/\mu$ we still obtain a valid skew
shape $(\lambda/\mu)^-$.  Suppose $\theta^-_\p$ has been defined for
$(\lambda/\mu)^-$.  If at most one of $(i-1,j)$ or $(i,j-1)$ is in
$(\lambda/\mu)^-$ then one can define $\theta_\p$ by making
$\theta_\p(i,j)$ either 1 or a very big value, letting
$\theta_\p(i',j') = \theta^-_\p(i',j')$ for other boxes $(i',j')$
(we may have to shift the values of $\theta^-_\p$ to be able to set
$\theta_\p(i,j) = 1$).

So assume that $(i-1,j), (i,j-1) \in (\lambda/\mu)^-$.  If
$O_\p((i-1,j) \lessdot (i,j)) = O_\p((i,j-1) \lessdot (i,j))$, then
$\theta_\p$ can be defined as in the previous case.  So assume
$O_\p((i-1,j) \lessdot (i,j)) = \overline{O_\p((i,j-1) \lessdot
(i,j))}$. If $(i-1,j-1) \notin (\lambda/\mu)^-$ then
$(\lambda/\mu)^-$ is disconnected.  In this case, we may pick
labelings $\theta^1_\p, \theta^2_\p$ for the two components $C_1,
C_2$ of $(\lambda/\mu)^-$ so that we can set $\theta_\p(C_1) =
\theta^1_\p(C_1)
> \theta_\p(i,j) > \theta_\p(C_2) = \theta^1_\p(C_2)$.

Finally, suppose that $(i-1,j-1) \in (\lambda/\mu)^-$ and assume
without loss of generality that $O_\p((i-1,j) \lessdot (i,j)) =
{\strict} = O_\p((i-1,j-1) \lessdot (i,j-1))$ and $O_\p((i,j-1)
\lessdot (i,j)) = {\weak} = O_\p((i-1,j-1) \lessdot (i-1,j))$ (we
have used the definition of $O_\p$).  Suppose $\theta^-_\p$ is
defined. Then $\theta^-_\p(i-1,j)>\theta^-_\p(i-1,j-1) >
\theta^-_\p(i,j-1)$. It suffices to define $\theta_\p(i,j)$ to be an
integer very close to $\theta^-_\p(i-1,j-1)$ and $\theta_\p(i',j') =
\theta^-_\p(i',j')$ for other boxes $(i',j')$, possibly shifting the
values so that $\theta^\p(i,j)$ can be inserted.

%There exists an extension of an $A$-labeled skew shape such that
%edge is weak if and only if its larger end is assigned a smaller
%number in the extension. Thus, corresponding $s^\p_{\lambda/\mu}$ is
%$L$-nonnegative.
%
%View the shape $\lambda/\mu$ with labeling $\p$ as a generalised
%poset. Reverse the inequalities on the strict edges of the Hasse
%diagram. We argue that the resulting inequalities are not in
%contradiction with each other. Any such contradiction would be of
%the form of a loop in Hasse diagram graph with each next vertex
%being larger than the previouse one, we call such loop {\it
%{tailed}}. Pick the vertex $(i,j)$ in this loop with maximal content
%$j-i$. Since $\p$ is alternating, one of the two edges of the loop
%ending in this vertex is reversed, while the other is not. Then this
%vertex is either smaller than both its neighbors in the loop, or
%larger - which means the loop is not tailed - contadiction.
%
%Now the statement of the Theorem follows from general theory of
%$P$-partitions, see Corollary 7.19.5, \cite{EC2}
\end{proof}

\begin{example}
\label{ex:lct6}
In Figure \ref{lct6} an edge-labeling $O_\p$ of the shape $\lambda = (2,2,1)$ is given. Here $p_{-1} = \weak$, $p_0 = \strict$, $p_1 = \weak$. The corresponding wave Schur function $s^\p_{\lambda}$ can be computed to be equal to $L_{(2,1,2)}+L_{(2,1,1,1)}+L_{(3,2)}+L_{(3,1,1)}+L_{(2,2,1)}$.  It is easy to check that this edge-labeling does come from a vertex labeling of the underlying poset.

\begin{figure}[h!]
\begin{center}
\input{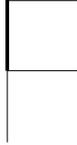}
\end{center}
\caption{An edge-labeling $O_\p$ of the shape $\lambda = (2,2,1)$.} \label{lct6}
\end{figure}

\end{example}

%\remind{complete the following remark}
\begin{remark}
Proposition~\ref{glp} implies a formula for
$s^\p_{\lambda/\mu}(1,q,q^2,q^3,\ldots)$ similar to that
of Proposition 7.19.11 in \cite{EC2}. Indeed, we can consider the descent set $D_{\p}(T)$ of a standard tableau $T$ with respect to $\theta_{\p}$. Then if we define a generalization of comajor index ${\rm comaj}_{\p}(T) = \sum_{i \in D_{\p}(T)} (n-i)$, we obtain the formula $$s^\p_{\lambda/\mu}(1,q,q^2,q^3,\ldots) = \frac{\sum_T q^{{\rm comaj}_{\p}(T)}}{(1-q)(1-q^2) \cdots (1-q^n)}.$$ However, it seems unlikely that an analog of hook content formula (see \cite[Theorem 15.3]{Sta71}) exists because the number of wave $\p$-tableaux filled with entries from $1$ to $n$ does not appear to factor nicely.  In Example \ref{ex:lct6} the number of wave $\p$-tableaux with entries from $1$ to $4$ is the prime number $23$.

\end{remark}

\begin{corollary}[Cell transfer for wave Schur functions]
\label{cor:cell} Let $\lambda/\mu$ and $\nu/\rho$ be two skew shapes
and $\p$ be any ``strict-weak'' assignment.  Then the difference
$s^\p_{\lambda/\mu \wedge \nu/\rho}\, s^\p_{\lambda/\mu
\vee\nu/\rho} - s^\p_{\lambda/\mu}\, s^\p_{\nu/\rho}$ is
$L$-positive.
\end{corollary}
\begin{proof}
Follows immediately from Theorem~\ref{thm:main} and
Proposition~\ref{glp}.
\end{proof}

In~\cite{LPP} it is shown that the difference in
Corollary~\ref{cor:cell} is in fact {\it Schur-positive} when $\p$
is the standard assignment.

\subsection{Jacobi-Trudi formula for wave Schur functions}

Let $\lambda = (\lambda_1,\lambda_2,\ldots,\lambda_l)$ and $\mu =
(\mu_1, \mu_2,\ldots,\mu_l)$ be two partitions satisfying $\mu
\subset \lambda$.  Now, for each pair $1 \leq i,j \leq l$ such that $\mu_j - j + 1 < \lambda_i - i$, define the set $$D_{ij}(\lambda,\mu) = \{\mu_j - j + 1 < a \leq \lambda_i - i \mid p_a = {\strict}\} - (\mu_j - j + 1).$$  Set
$\alpha_{ij}(\lambda,\mu) = \C(D_{ij}(\lambda,\mu))$ to be the
corresponding composition of $\lambda_i - \mu_j - i + j$. If $\mu_j - j + 1 = \lambda_i - i$, set $\alpha_{ij}(\lambda,\mu) = (1)$. If $\mu_j - j = \lambda_i - i$ set $\alpha_{ij}(\lambda,\mu) = (0)$. Finally, if $\mu_j - j > \lambda_i - i$ set $\alpha_{ij}(\lambda,\mu) = \emptyset$. Let $L_{(0)}=1$, $L_{\emptyset} = 0$.

%If Start from the leftmost cell in row $i$ of $(\lambda/\mu)_\p$ and
%go along one of the possible shortest paths without left or down
%steps to the diagonal of the rightmost cell in the $j$-th row of
%$(\lambda/\mu)_\p$. If the content of the rightmost cell in the
%$j$-th row is less then content of the leftmost cell in row $i$
%minus $1$, then the entry is $0$. At last, if the content of the
%rightmost cell in the $j$-th row is is equal to the content of the
%leftmost cell in row $i$ minus $1$, the entry of the matrix is $1$.

\begin{theorem}[Jacobi-Trudi expansion for wave Schur functions]
\label{jtl} Let $\lambda/\mu$ be a skew shape.  Then
$$
s^\p_{\lambda/\mu} = {\rm det}(L_{\alpha_{ij}(\lambda,\mu)})_{i,j =
1}^n
$$
where $n$ is the number of rows in $\lambda$.
\end{theorem}

\begin{example}
Let $\lambda = (7,6,6,4)$, $\mu = (2,2,1,0)$. Then for $\p$ given by $$\ldots p_{-3}, p_{-2}, p_{-1}, p_{0}, p_{1}, p_{2}, p_{3},
p_{4}, p_{5}, p_{6}, p_{7}, p_{8}, \ldots =$$ $$\ldots \strict,
\strict, \strict, \weak, \strict, \weak, \weak, \strict, \strict,
\weak,\strict,\weak \ldots$$ we get the shape in Figure~\ref{fig:gens}, and

\begin{eqnarray*}
s^\p_{\lambda/\mu} =
\left|
\begin{array}{cccc}
L_{(2,1,2)} &  L_{(3,1,2)} &  L_{(2,3,1,2)}  &  L_{(1,1,2,3,1,2)} \\
L_{(2,1)} & L_{(3,1)} &  L_{(2,3,1)}  &  L_{(1,1,2,3,1)} \\
L_{(2)} & L_{(3)} &  L_{(2,3)}  &  L_{(1,1,2,3)} \\
0 & 1 &  L_{(2)}  &  L_{(1,1,2)}
\end{array}
\right|
.
\end{eqnarray*}
%
%Here to obtain third entry of the first row of the matrix we start
%reading the third row of diagram at point $(3,2)$, which gives us
%weak, strict, weak, weak. Then we read two vertical edges from
%$(3,6)$ to $(1,6)$, and since we reverse the labels we get strict,
%strict. Finally, we get one more weak edge from the first row,
%ending in $(1,7)$. Combining, we get $L_{(2,3,1,2)}$ entry in the
%matrix. Another example: the leftmost point in second row has
%content $1$, while the rightmost point of forth row has content $0$,
%exactly one less. Therefore the corresponding entry is $1$.
\end{example}

\begin{figure}
\begin{center}
\input{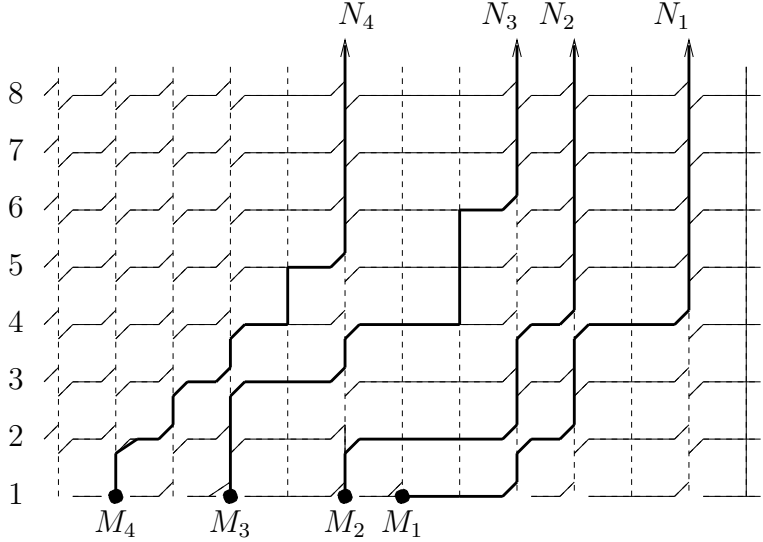}
\end{center}
 \caption{A family of paths on altered grid corresponding to $\p$-tableau on Figure \ref{lct8}.}
\label{fig:qsct}
\end{figure}

\begin{proof}[Proof of Theorem~\ref{jtl}]
Let us construct an oriented network $N_\p$, which depends on the
choice of $\p$.  Namely, we begin with the square grid built on the
points in the upper half plane, with row 1 being the bottom row, and
orient all edges to the right or upwards.  Then we alter the all the
crossings in each column $C_i$ such that $p_i = \strict$ as shown in Figure \ref{lct7}. Namely, we arrange these crossings so that it is impossible to move
from left to right through them, but other directions that were
possible before are still possible (see Figure~\ref{fig:qsct}).  We
assign to each edge in row $i$ weight $x_i$, and every other edge
weight $1$. Now mark the points $M_k$ with coordinates $(\mu_k - k +
1, 1)$ on our grid.  Mark exit directions $N_k$ in the columns
numbered $\lambda_k - k + 1$.

\begin{figure}[h!]
\begin{center}
\input{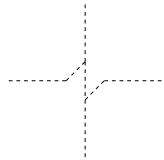}
\end{center}
\caption{A local picture of an altered crossing.} \label{lct7}
\end{figure}

Now we apply the Gessel-Viennot method to this path network; see for
example~\cite[Chapter 7]{EC2} for the application of this method in
the case of Schur functions.  For each pair $1 \leq i,j \leq n$ the
weight generating function of the paths from $M_i$ to $N_j$ is equal
to $L_{\alpha_{ij}}(\lambda,\mu)$. Thus the determinant ${\rm
det}(L_{\alpha_{ij}(\lambda,\mu)})_{i,j = 1}^n$ is equal to the
weight generating function of families of non-crossing paths
starting at the $M_i$-s and ending in the columns $N_i$. These
families of non-crossing paths are in (a weight-preserving)
bijection with wave $\p$-tableau of shape $\lambda/\mu$.
\end{proof}

%\remind{Figure needs a caption.  Draw a picture of what happens
%locally at a single vertex, add notation.  Also need arrows on
%picture.  Draw the tableau corresponding to the above family of
%paths.}

\begin{figure}[h!]
\begin{center}
\input{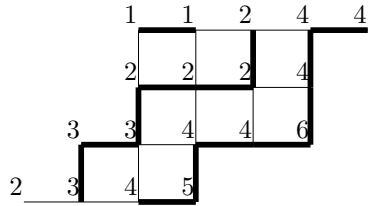}
\end{center}
\caption{A wave $\p$-tableau of the shape $\lambda/\mu = (7,6,6,4)/(2,2,1,0)$ with the edge labeling $O_{\p}$ as in Figure \ref{fig:gens}.} \label{lct8}
\end{figure}

\begin{remark}
(i) We have $\omega(s^\p_{\lambda/\mu}) = s^{\bar \p}_{\lambda/\mu}$
where $\bar \p = (\ldots, \overline{p_{-2}}, \overline{p_{-1}},
\overline{p_0}, \overline{p_1}, \ldots)$.

(ii) Let us denote $\widetilde {\lambda/\mu}$ the rotated on $180$ degrees $\p$-tableau $\lambda/\mu$ with $\tilde p_i = p_{-i}$. Then $\nu(s^\p_{\lambda/\mu}) = s^{\tilde \p}_{\widetilde {\lambda/\mu}}$.
\end{remark}

The following theorem, combined with Proposition~\ref{glp}, answers
Question~\ref{q:main} for the case that $Q$ and $R$ are convex
subsets of a chain.
\begin{theorem}
\label{thm:diff} The differences $L_{\alpha \wedge_m \beta}\,
L_{\alpha \vee_m \beta} - L_\alpha \, L_\beta$ of
Theorem~\ref{thm:lp} are equal to wave Schur functions.
\end{theorem}

\begin{figure}
\begin{center}
\input{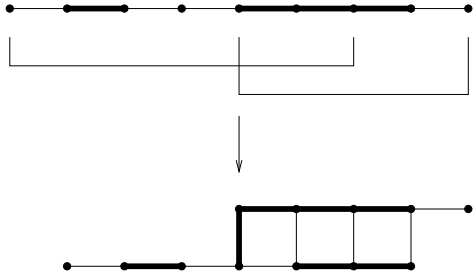}
\end{center}
\caption{The skew shape corresponding to the difference $L_{(1,2,1,3)}L_{(4,1)}-L_{(3)} L_{(1,2,1,4,1)}$.}
\label{fig:JT}
\end{figure}

\begin{proof}
We may suppose that $m \geq 1$ for otherwise the difference is 0.
Pick a sequence $\p = \p^\alpha$ such that $p_i = \strict$ if and
only if $i \in D(\alpha)$ (this determines $p_1, p_2, \ldots,
p_{|\alpha|-1}$). Then set $\lambda = (|\alpha|,m + |\beta|)$ and $\mu =
(m -1, 1)$.

We can compute that
\begin{align*}
L_{\alpha_{11}(\lambda,\mu)} &= L_{\alpha \vee_m \beta} &
L_{\alpha_{12}(\lambda,\mu)} &= L_{\alpha} \\
L_{\alpha_{21}(\lambda,\mu)} &= L_{\alpha \wedge_m \beta} &
L_{\alpha_{22}(\lambda,\mu)} &= L_{\beta}.
\end{align*}

Theorem~\ref{jtl} tells us that $s_{\lambda/\mu}^\p$ is exactly
$\det\left(
\begin{array}{cc}
L_{\alpha \vee_m \beta} &  L_{\alpha}\\
L_{\beta} & L_{\alpha \wedge_m \beta}
\end{array}
\right)$.

%Form a skew shape from $\alpha \wedge_m \beta$ and $\alpha \vee
%\beta$ as follows: place $\alpha \wedge \beta$ one row under $\alpha
%\vee \beta$ so that the vertecies coming from $\beta$ have the same
%content as their pre-images. Consider the union of this shifted
%$\alpha \wedge \beta$ and original $\alpha \vee \beta$, and label
%vertical edges between them so that we get an $A$-labeling.
\end{proof}
We illustrate the choice of $\p^\alpha$, $\lambda$ and $\mu$ of
Theorem~\ref{thm:diff} in Figure~\ref{fig:JT}. Here $\p^\alpha$ is such that $\ldots p_1, p_2, p_3, p_4, p_5, p_6, p_7, p_8, \ldots =$ $$\ldots, \strict, \weak, \strict, \strict, \weak, \weak, \weak, \strict,$$ $\alpha = (1,2,1,4,1)$, $\beta = (3)$, $m=5$. Then $\lambda = (9,8)$, $\mu = (4,1)$ and the corresponding $s^\p_{\lambda/\mu}$ is equal to $L_{(1,2,1,3)}L_{(4,1)}-L_{(3)} L_{(1,2,1,4,1)}$.

%\remind{Say what
%$\p^\alpha$, $\lambda$ and $\mu$, $\alpha$, $\beta$, $m$, ... are}

\section{More general labeled posets}
\label{sec:gen}
Our point of view so far has been that the $P$-partition generating
functions $K_{P,\theta}$ are ``skew'' analogues of the fundamental
quasi-symmetric functions $L_\alpha$, just as skew Schur functions
are skew versions of the usual Schur functions.  From this point of
view, the two key properties that the generating functions
$K_{P,\theta}$ possess are (a) they lie in $\QSym$, and (b) they are
$L$-positive.

\medskip
\subsection{$\T$-labeled posets}

In~\cite{LP}, we defined more general {\it $\T$-labeled posets} for
which Theorem~\ref{thm:celltransfer} also holds.

Let $\T$ denote the set of all weakly increasing functions $f: \P
\rightarrow \Z \cup \{\infty\}$.  A {\it $\T$-labeling} $O$ of a
finite poset $P$ is a map $O : \{(s,t)\in P^2 \mid s \gtrdot t \}
\rightarrow \T$ labeling each edge $(s,t)$ of the Hasse diagram by a
weakly increasing function $O(s,t):\P \rightarrow \Z \cup
\{\infty\}$.  A {\it $\T$-labeled poset} is an an ordered pair
$(P,O)$ where $P$ is a poset, and $O$ is a $\T$-labeling of $P$.  It
is clear how to take directed sums of $\T$-labeled posets, or to
take convex subsets of $\T$-labeled posets.

A {\it $(P, O)$-tableau} is a map $\sigma: P \to \mathbb P$ such
that for each covering relation $s \lessdot t$ in $P$ we have
$\sigma(s) \leq O(s,t)(\sigma(t))$.  Note that ``tableau'' here is
used in the same sense as ``partition'' was in
Section~\ref{sec:Ppart}.  The weight generating function of all
$(P,O)$-tableaux is denoted $K_{P,O}$.

\begin{problem}
\label{p:qsym} For which $\T$-labeled posets is the generating
function $K_{P,O}$ quasi-symmetric?
\end{problem}
We call $(P,O)$ {\it quasi-symmetric} if $K_{P,O} \in \QSym$.  There
is a large class of quasi-symmetric $\T$-labeled posets, containing
all those induced from the form $(P,\theta)$.  Define $f^{\weak}(x)
= x$ and $f^{\strict}(x) = x-1$.  Then $K_{P,O} \in \QSym$ if $O(s
\lessdot t) \in \{f^{\weak}, f^{\strict}\}$ for every covering
relation $s \lessdot t$. Following terminology of
McNamara~\cite{McN}, such strict-weak edge labeled posets are called
{\it oriented}.

It is unclear how to obtain more solutions to Problem~\ref{p:qsym}.
However, we can show, using a factorization result we prove later,
 that that the answer to Problem~\ref{p:qsym}is
compatible with taking disjoint unions and connected components.

\begin{proposition}
\label{prop:connect} If $(P,O_P)$ and $(Q,O_Q)$ are quasi-symmetric
then so is $(P \oplus Q, O_P \oplus O_Q)$.  If $(P,O)$ is
quasi-symmetric then each connected component $(P_i,O|_{P_i})$ of
$(P,O)$ is also quasi-symmetric.
\end{proposition}
\begin{proof}
Since $K_{P \oplus Q, O_P \oplus O_Q} = K_{P,O}\,K_{Q,O}$, the first
statement holds because $\QSym$ is a ring.  The second statement
holds by Theorem~\ref{thm:fact}.
\end{proof}

In particular, Proposition~\ref{prop:connect} says that it is not
possible to obtain a quasi-symmetric $\T$-labeled poset
``accidentally'' by taking disjoint sums.

\subsection{Oriented posets}
Let us say that an orientation $O$ of a poset $P$ arises from a
labeling $\theta$ if $O(s \lessdot t) = f^{\strict}$ exactly when
$(s \lessdot t)$ is a descent of $\theta$.  Clearly in this case we
have $K_{P,O} = K_{P,\theta}$.  Not every orientation arises from a
labeling, as shown in~\cite[Example 2.7]{McN}.  It is also possible
to find both oriented posets such that the generating function 
$K_{P,O}$ is $L$-positive and
oriented posets such that $K_{P,O}$ is not $L$-positive;
see~\cite[Remark~5.9]{McN}.

\begin{problem}
\label{p:lpos} For which oriented posets $(P,O)$ is the generating function
$K_{P,O}$ $L$-positive?
\end{problem}

It is unclear to us whether an analogue of Theorem~\ref{thm:main}
should hold for more general oriented posets $(P,O)$, such as the
ones which are solutions to Problem~\ref{p:lpos}.  This would
potentially expand our notion of ``skew'' fundamental
quasi-symmetric functions beyond just $K_{P,\theta}$.

\section{Algebraic properties of $\QSym$} 
\label{sec:fact}
We
prove in this section some algebraic results concerning $\QSym$ used
earlier.

\subsection{A factorization property of quasi-symmetric functions}
\label{fac} Denote by $K = \Z[[x_1,x_2,x_3,\ldots]]$ the ring of
formal power series in infinitely many variables with bounded
degree.  Clearly the units in $K$ or in $K^{(n)}$ are $1$ and $-1$.
In this subsection, we prove the following property of $\QSym$.

\begin{theorem}
\label{thm:fact} Suppose $f \in \QSym$ and $f = \prod_i f_i$ is a
factorization of $f$ into irreducibles in $K$.  Then $f_i \in \QSym$
for each $i$.
\end{theorem}

Now let $\a = (1 \leq a_1 < a_2 < \ldots <)$ be an increasing
sequence of positive integers and let $A$ denote the set of such
sequences. Define the algebra homomorphism $A_\a: K \to K$ by
$$
A_\a f := f(x_\a) := f(0,\ldots,0, x_1,0,\ldots,0,x_2,0,\ldots)
$$
where $x_i$ is placed in the $a_i$-th position.  For a sequence
$\a$, we shall also write $a(i) = a_i$ in function notation.  Thus
$a: {\mathbb N} \to {\mathbb N}$ is a strictly increasing function.

As an example, take $a = (2,3,4,\dotsm)$, $b = (1,3,5,\dotsm)$. Then we have
$A_a f = f(0,x_1,x_2,\ldots)$ and $A_b \circ A_a f = A_b (A_a f) =
f(0,x_1,0,x_2,0,x_3,\dotsm)$.

The following lemma is essentially the definition.
\begin{lemma}
\label{lem:Qsym} An element $f \in K$ is quasi-symmetric if and only
if $f(x_\a) = f$ for each $\a \in A$.
\end{lemma}

Let $k \geq 1$ be an integer.  Define $\a^{(k)}$ by $$a^{(k)}(i) =
\begin{cases} i & \mbox{if $i < k$,} \\ i+1 & \mbox{if $i \geq
k$.}\end{cases}$$

\begin{lemma}
\label{lem:Qsymn} Suppose $f \in K$ has degree $n$.  Then $f$ is
quasi-symmetric if and only if $f(x_\a) = f$ for the sequences
$\a^{(k)}$ for $1 \leq k \leq n$.
\end{lemma}
\begin{proof}
The only if direction is clear.  Assume that $f(x_\a) = f$ for each
$\a^{(k)}$ for $1 \leq k \leq n$.  To show that the coefficients of
$x_1^{c_1} \cdots x_n^{c_n}$ and $x_{b_1}^{c_1} \cdots
x_{b_n}^{c_n}$ in $f$ are the same we use (the coefficient of
$x_1^{c_1} \cdots x_n^{c_n}$  in the equality)
$$
A_{\a^{(n)}}^{b_n-b_{n-1}-1} \cdots
A_{\a^{(2)}}^{b_2-b_1-1}A_{\a^{(1)}}^{b_1-1}
 f = f.
$$
\end{proof}

The following lemma is a simple calculation.
\begin{lemma}
\label{lem:combine} We have $A_\b \circ A_\a = A_\c$ where $c(i) =
a(b(i))$.
\end{lemma}

\begin{lemma}
\label{lem:finite} Let $f \in K$.  Suppose $f$ has finite order with
respect to $A_\a$ for every $\a \in A$.  Then there exists $\b \in
A$ so that $A_\b f \in \QSym$.
\end{lemma}
\begin{proof}
Given $f \in K$ invariant under $A_{\a^{(k)}}$ for $1 \leq k \leq
t$, with $t$ possibly $0$, we will produce an $f' = A_\b f$
invariant under $A_{\a^{(k)}}$ for $1 \leq k \leq t+1$.  Using
Lemma~\ref{lem:Qsymn} and the fact that $f$ has bounded degree this
is sufficient.

So let $f$ be invariant under $A_{\a^{(k)}}$ for $1 \leq k \leq t$.
By assumption $A_{\a^{(t+1)}}$ has finite order $d$ on $f$.  Define
$\b \in A$ by $b(i) = 1 + (i-1)d$ and let $f' = A_\b f$.  We claim
that $f'$ is invariant under  $A_{\a^{(k)}}$ for $1 \leq k \leq
t+1$.  We have
$$
b(a^{(k)}(i)) = \begin{cases} 1+(i-1)d & \mbox{if $i < k$,} \\
1+id & \mbox{if $i \geq k$.}\end{cases}
$$
In the following we will repeatedly use Lemma~\ref{lem:combine}.

Define $\b^{(j)} \in A$ for $1 \leq j < k$ by
$$
b^{(j)}(i) = \begin{cases} i & \mbox{if $i < j$,} \\ j + (i-j)d &
\mbox{$ i \geq j$.} \end{cases}
$$
Note that $A_{\b^{(j)}} \circ (A_{\a^{(j)}})^{d-1} =
A_{\b^{(j-1)}}$. Similarly define $\c^{(j)} \in A$ for $1\leq j < k$
by
$$
c^{(j)}(i) = \begin{cases} i & \mbox{if $i \leq j$,} \\ j + (i-j)d &
\mbox{$ j < i < k$,} \\ j+ (i-j+1)d & \mbox{$i \geq k$.}\end{cases}
$$
Note that $A_{\c^{(j)}} \circ (A_{\a^{(j)}})^{d-1} =
A_{\c^{(j-1)}}$. We also have the three equalities \begin{align*}
A_{\c^{(k-1)}} &= A_{\b^{(k-1)}}
\circ (A_{\a^{(k)}})^d, &
A_{\b^{(1)}} &= A_\b, &
 \text{and} \;\;A_{\c^{(1)}} &= A_{\a^{(k)}} \circ A_\b.
\end{align*}

Finally using our assumptions and $ 1 \leq k \leq t+1$, we have
$$
A_\b f = A_{\b^{(1)}} f = \cdots =  A_{\b^{(k-1)}}f =
A_{\c^{(k-1)}} f = \cdots =  A_{\c^{(1)}}f = A_{\a^{(k)}} \circ A_\b
f.$$
\end{proof}

\begin{proof}[Proof of Theorem~\ref{thm:fact}]
Let $\a \in A$.  Applying $A_\a$ to $f = \prod_i f_i$ and using
Lemma~\ref{lem:Qsym}, we have $f = \prod_i A_\a f_i$.  By
Lemma~\ref{lem:KUFD} below, each $A_\a f_i$ must be equal to $\pm f_j$. In
other words, $A_\a$ has finite order on each $f_i$ and application
of $A_\a$ to $f_i$ produces (up to sign) another $f_j$.  Using
Lemma~\ref{lem:finite}, we see that $f_j$ must lie in $\QSym$ for
some $j$.  Now divide both sides by $f_j$ and proceed by induction.
\end{proof}

\begin{corollary}
\label{cor:UFD} $\QSym$ is a unique factorization domain.
\end{corollary}

Corollary~\ref{cor:UFD} also follows from work of
Hazewinkel~\cite{Haz}, who shows that $\QSym$ is a polynomial ring.

\begin{proof}
If $f \in \QSym$ then two irreducible factorizations of $f$ in
$\QSym$ will also be irreducible factorizations in $K$, by
Theorem~\ref{thm:fact}.  The theorem follows from
Lemma~\ref{lem:KUFD}, proven below.
\end{proof}

%The following result is probably well known, but we have been unable to find a suitable reference.
\begin{lemma}
\label{lem:KUFD} The ring $K$ is a unique factorization domain.
\end{lemma}
\begin{proof}
We start by recalling the well known fact that the polynomial rings $K^{(n)} =
\Z[x_1,x_2,\ldots,x_n]$ are unique factorization domains. An element
of $f(x_1,x_2,\ldots) \in K$ is determined by its images
$$f^{(n)} = f(x_1,x_2,\ldots,x_n,0,0,\ldots) \in K^{(n)}.$$  We may
write $f = (f^{(n)})$ for a compatible sequence of $f^{(n)} \in
K^{(n)}$.

We first claim that $f$ is irreducible if and only if there exists
$N > 0$ such that $f^{(n)}$ is irreducible for all $n > N$. Let $f =
\prod f_i$ be a decomposition of $f$ into irreducibles.  Then there
exists $M > 0$ so that $\deg(f_i^{(n)}) = \deg(f_i)$ for all $i$ and
$n > M$.  Thus $f^{(n)} = \prod f_i^{(n)}$ is reducible for $n > M$
if $f$ is.  Conversely, suppose that $f^{(n)}$ is reducible for
infinitely many values of $n$.  If $n > M$ and $f^{(n)}$ is
reducible then $f^{(m)}$ is also reducible for $n > m > M$.  Thus we
may assume $f^{(n)}$ is reducible for all $n > N$ for some $N > M$.
Restriction of $f^{(n)}$ to $f^{(m)}$ for $n > m > N$ will not
change the degree of any of the factors. Thus the factorizations of
$f^{(n)}$ are compatible for each $n > N$.  For sufficiently large
$n \gg N$, the number $k$ of irreducible factors of $f^{(n)}$ will
be constant and greater than 1.  Ordering the factorizations
$\prod_{i=1}^k f^{(n)}_i$ compatibly, we conclude that $f =
\prod_{i=1}^k f_i$ where $f_i = (f^{(n)}_i)$ is reducible.

Now suppose that $f = \prod_i f_i = \prod_j g_j$ are two
factorizations of $f$ into irreducibles.  By our claim, there exists
some huge $N$ so that $\prod_i f_i^{(n)} = \prod_j g_j^{(n)}$ are
factorizations of $f^{(n)}$ into irreducibles in $K^{(n)}$, for each
$n > N$. Since $K^{(n)}$ is a $UFD$, these factorizations are the
same up to permutation and sign: $g^{(n)}_i = \epsilon_i
f^{(n)}_{\sigma(i)}$. If $N$ is chosen large enough the same
permutation $\sigma$ and signs $\epsilon_i$ will work for all $n >
N$.  This shows that $g_i = \epsilon_i f_{\sigma(i)}$.
\end{proof}

\begin{remark}
(i) Note that Corollary~\ref{cor:UFD} is not true in finitely many
variables. For example, in two variables $x_1$ and $x_2$ we have
$(x_1^2 x_2) (x_1x_2^2) = (x_1x_2)^3$.

(ii) It seems interesting to ask whether the $r$-quasi-symmetric
functions defined by Hivert~\cite{Hiv} also form a unique factorization
domain. The $m$-quasi-invariants~\cite{EG} occurring in
representation theory do not in general form unique factorization
domains.
\end{remark}

\subsection{Irreducibility of fundamental quasi-symmetric functions}
In this section, we show that the fundamental quasi-symmetric functions
$\{L_{\alpha}\}$ and the monomial quasi-symmetric functions
$\{M_{\alpha}\}$ are irreducible in $\QSym$ and in $K$.

Let $\alpha = (\alpha_1,\alpha_2,\ldots,\alpha_k)$ and $\beta =
(\beta_1,\beta_2,\ldots,\beta_l)$ be two compositions.  Define the
{\it lexicographic} order on compositions by $\alpha
> \beta$ if and only if for some $i$ we have $\alpha_j = \beta_j$
for $1 \leq j \leq i-1$ and $\alpha_i > \beta_i$.  Using this order,
we obtain lexicographic orders on monomials $\{x^\alpha\}$, monomial
quasi-symmetric functions $\{M_\alpha\}$ and fundamental
quasi-symmetric functions $\{L_\alpha\}$.  Note that the
lexicographically maximal monomial in $M_\alpha$ or $L_\alpha$ is
$x^\alpha$.

In the following proofs we say that a quasi-symmetric function $f$
contains a term $L_\alpha$ (and similarly for $M_\alpha$) if the
coefficient of $L_\alpha$ is non-zero when $f$ is written in the
basis of fundamental quasi-symmetric functions.  The following lemma
is immediate from the definitions.

\begin{lemma}
\label{lem:lex} The lexicographically maximal monomial in the
product $f\,g$ of two quasi-symmetric functions $f$ and $g$ is the
product of the lexicographically maximal monomials in $f$ and $g$.
\end{lemma}

\begin{proposition}
\label{prop:mirred}
The monomial quasi-symmetric function $M_\alpha$ is irreducible in
$\QSym$ and in $K = \Z[[x_1,x_2,\ldots]]$.
\end{proposition}

\begin{proof}
We proceed by induction on the size $n = \alpha_1+ \alpha_2 + \cdots
+ \alpha_k$.  For $n=1$ the statement is obvious.

Assume now that $M_{\alpha} = f \, g$ is not irreducible.  Note
first that $f$ and $g$ must be homogeneous.  Otherwise, the
homogeneous components of maximal and minimal degree in the product
would not cancel out, and thus we would never get the homogeneous
function $M_{\alpha}$. Also note that according to
Theorem~\ref{thm:fact} both $f$ and $g$ must be quasi-symmetric.   First we suppose that $k = 1$.

Now, take the specialization $x_i = q^{i-1}$. It is known (\cite[Proposition 7.19.10]{EC2}) that under this specialization we have $L_{\alpha}(1,q,q^2,\ldots) = \frac{q^{e(\alpha)}}{(1-q)(1-q^2) \cdots (1-q^{|\alpha|})}$, where $e(\alpha)$ is the ``comajor'' statistic.  That means that if ${\rm deg}(f) = p$, ${\rm deg}(g) = \alpha_1 - p$ and $0 < p < \alpha_1$, $fg$ will never have a pole at primitive $\alpha_1$-th root of unity. On the other hand $M_{\alpha}(1, q, \ldots) = \frac{1}{1-q^{\alpha_1}}$, and thus $M_{\alpha}$ has a pole at a primitive $\alpha_1$-th root of unity, which is a contradiction.

%We
%may thus assume that $k > 1$, for otherwise we may apply the
%involution $\omega$ of Proposition~\ref{prop:inv} to obtain a
%factorization of $L_{\bar \alpha}$.

Now suppose that $k \neq 1$.  We write each of the participating functions as polynomials in
$x_1$:
\begin{align*}
M_{\alpha} &= x_1^{\alpha_1} \tilde M_1 + \cdots, &f &= x_1^r \tilde
f_1 + \cdots,  &g &= x_1^{\alpha_1-r} \tilde g_1 + \cdots.
\end{align*}
Here the leading term is the one with the highest power of $x_1$,
and the notation $\tilde f$ denotes a power series
$f(x_2,x_3,\ldots)$ quasi-symmetric in the variables $x_2, x_3,
\dotsm$.

Note that $\tilde M_1 = M_{(\alpha_2, \dotsm, \alpha_k)}(x_2,x_3,
\dotsm)$ is the monomial quasi-symmetric function corresponding
to the composition obtained from $\alpha$ by removing the first part.
Since $M_{\alpha} = f \, g$, we must have $M_1 = f_1 \, g_1$. By
induction one of $f_1$ or $g_1$ is equal to a unit, $\pm 1$. Without
loss of generality we can assume $f_1 = 1$.  Thus the monomial
quasi-symmetric function $M_{(r)}$ occurs in $f$.  By
Lemma~\ref{lem:lex} above we conclude that the lexicographically
maximal monomial quasi-symmetric function in $g$ is $M_{(\alpha_1-r,
\alpha_2, \dotsm, \alpha_k)}$.

%In particular this monomial has non-zero coefficient in $Q$.

Now apply the involution $\nu$ of Proposition~\ref{prop:inv} to the
equality $M_{\alpha} = f\, g$ to obtain $M_{\alpha^*} = \nu(f)
\nu(g)$.  By Proposition~\ref{prop:inv}, $\nu(M_{(r)}) = M_{(r)}$
and so the monomial $M_{(r)}$ is still the lexicographically maximal
monomial in $\nu(f)$. Similarly, the monomial symmetric function
$M_{(\alpha_k, \dotsm, \alpha_2, \alpha_1-r)}$ occurs in $\nu(g)$
with non-zero coefficient.  Since $k \neq 1$, the lexicographically
maximal monomial in the product $\nu(f)\nu(g)$ is at least as large
as $M_{(\alpha_k+r, \dotsm, \alpha_1-r)}$. This however is lexicographically
larger than $M_{\alpha^*} = M_{(\alpha_k, \dotsm, \alpha_1)}$ unless $r
= 0$.

We conclude that $f = 1$ and that $M_\alpha$ is irreducible.
\end{proof}

\begin{proposition}
\label{prop:irred} The fundamental quasi-symmetric function
$L_\alpha$ is irreducible in both $\QSym$ and in $K =
\Z[[x_1,x_2,x_3,\ldots]]$.
\end{proposition}

\begin{proof}
The trick used for the case $k = 1$ in the proof of Proposition~\ref{prop:mirred} also works here.

%The only claim which requires a new argument is irreducability in case $k=1$. Assume in this case $M_{\alpha} = \sum_i x_i^
%{\alpha_1}$ is equal to the product $fg$. Then by Theorem \ref{thm:fact} $f$ and $g$ are quasisymmetric. Both $f$ and $g$ can be
%assumed to be homogenous. Indeed, $\mathbb C[[x_1, x_2, \ldots]]$ is an integral domain. Therefore taking highest and lowest
% degree homogenous components of $f$ and $g$ we see that $fg$ would have more than one non-zero homogenouse component if %either of $f$ and $g$ did. 

%Now, take the specialization $x_i = q^{i-1}$. It is known (\cite{EC2}, Proposition 7.19.10) that under such specialization $L_{\alpha} = 
%\frac{q^{e(\alpha)}}{(1-q)(1-q^2) \cdots (1-q^{|\alpha|})}$, where $e(\alpha)$ is a certain statistics. That means that if $deg(f) = p$, 
%$deg(g) = \alpha_1 - p$ and $0 < p < \alpha_1$, $fg$ will never have a pole at primitive $\alpha_1$-th root of unity. On the other h
%and $M_{\alpha}(1, q, \ldots) = \frac{1}{1-q^{\alpha_1}}$, and thus $M_{\alpha}$ has a pole at primitive $\alpha_1$-th root of unity - 
%contradiction.
\end{proof}


\begin{thebibliography}{a-a-a-a-a-a-a}

%\bibitem[BS]{BS} {\sc N.~Bergeron and F.~Sottile:}
%Skew Schubert functions and the Pieri formula for flag manifolds,
%{\sl Trans. Amer. Math. Soc.}, \textbf{354} No.2, (2002), 651-673.

\bibitem[EG]{EG} {\sc P.~Etingof and V.~Ginzburg:}
On $m$-quasi-invariants of a Coxeter group,  {\sl Mosc. Math. J.}
\textbf{2} (2002), 555--566.

\bibitem[Haz]{Haz} {\sc M.~Hazewinkel:}
The algebra of quasisymmetric functions is free over the integers,
{\sl Advances in Mathematics}, {\bf {164}}, 2001, 283-300.

\bibitem[Hiv]{Hiv} {\sc F.~Hivert:} Local action of the symmetric
group and generalizations of quasi-symmetric functions, {\sl. Proc.
of FPSAC}, Vancouver, 2004.

\bibitem[LP]{LP} {\sc T.~Lam and P.~Pylyavskyy:} Cell transfer and
monomial positivity, arXiv: {\tt math.CO/ 0505273}.

\bibitem[LPP]{LPP} {\sc T.~Lam, A.~Postnikov and P.~Pylyavskyy:}
Schur positivity and Schur log-concavity, preprint; {\tt
math.CO/0502446}.

\bibitem[LPP2]{LPP2} {\sc T.~Lam, A.~Postnikov and P.~Pylyavskyy:}
Some Schur positivity conjectures, in preparation.

\bibitem[McN]{McN} {\sc P.~McNamara:} Cylindric Skew Schur
Functions, {\sl Adv. Math.}, to appear; {\tt math.CO/0410301}.

%\bibitem[Pos]{Pos} {\sc A.~Postnikov:} Affine approach to quantum
%Schubert calculus, {\sl Duke Math. J.}, to appear; {\tt
%math.CO/0205165}.

\bibitem[Sta99]{EC2}{\sc R.~Stanley:} {\sl Enumerative
Combinatorics, Vol 2}, Cambridge, 1999.

\bibitem[Sta71]{Sta71} {\sc R.~Stanley:} Theory and applications of plane partitions: Part II, {\sl Stud. Appl. Math.} no. 50 (1971), 259-279.

\bibitem[Sta72]{Sta} {\sc R.~Stanley:} Ordered structures and partitions,
{\sl Memoirs Amer. Math. Soc.,} no. 119 (1972).

\end{thebibliography}
\end{document}